\documentclass[11pt]{article}
\usepackage[top=0.8in, bottom=0.8in, left=1.25in, right=1.25in]{geometry}

\usepackage{amsmath,amssymb,amsthm,xcolor,enumerate}
\usepackage{graphicx}
\usepackage[unicode,breaklinks=true,colorlinks=true,citecolor = {magenta}]{hyperref}
\graphicspath{ {./photopng/} }
\usepackage{float}

\newcommand{\dis}{\displaystyle}
\newcommand{\donothing}[1]{{}}

\numberwithin{equation}{section}
\newtheorem{theorem}{Theorem}[section]

\newtheorem{lemma}[theorem]{Lemma}
\newtheorem{proposition}[theorem]{Proposition}

\theoremstyle{remark}
\newtheorem{remark}[theorem]{Remark}
\newtheorem{example}[theorem]{Example}
\newcommand{\bke}[1]{\left ( #1 \right )}

\newcommand{\bket}[1]{\left \{ #1 \right \}}
\newcommand{\norm}[1]{ \| #1  \|}

\newcommand{\abs}[1]{\left | #1 \right |}
\newcommand{\mat}[1]{\begin{bmatrix} #1 \end{bmatrix}}

\newcommand\al{\alpha}

\newcommand\ga{\gamma}
\newcommand\de{\delta}

\newcommand\e {\varepsilon}

\newcommand\la{\lambda}
\newcommand\si{\sigma}
\newcommand\om{\omega}
\newcommand\Ga{\Gamma}
\newcommand\De{\Delta}

\newcommand\Si{\Sigma}
\newcommand\Om{\Omega}

\newcommand{\R}{\mathbb{R}}
\newcommand{\CC}{\mathbb{C}}

\newcommand{\N}{\mathbb{N}}

\newcommand{\cO}{\mathcal{O}}

\newcommand{\sech} {{\mathrm{sech}}}

\newcommand{\sgn}  {\mathop{\mathrm{sgn}}}

\newcommand{\pd}{\partial}

\newcommand{\ch}{\mathrm{ch}}
\newcommand{\sh}{\mathrm{sh}}

\newcommand{\I}{\infty}

\newcommand{\EQ}[1]{\begin{equation} #1 \end{equation}}
\newcommand{\EQS}[1]{\begin{equation}\begin{split} #1 \end{split}\end{equation}}
\newcommand{\EQN}[1]{\begin{equation*}\begin{split} #1 \end{split}\end{equation*}}
\newcommand{\EN}[1]{\begin{enumerate} #1 \end{enumerate}}
\newcommand{\no}{\textup{no}}
\newcommand{\ex}{\textup{ex}}

\title{Existence and stability of standing waves for one dimensional
  NLS with triple power nonlinearities}

\author{Fei Liu\thanks{Department of Applied Mathematics, University of Waterloo, Waterloo, ON N2L 3G1, Canada.\newline
\texttt{feiliu0625@gmail.com}}
\and
Tai-Peng Tsai\thanks{Department of Mathematics, University of British Columbia, Vancouver, BC V6T 1Z2, Canada.\newline ttsai@math.ubc.ca.
}
\and
Ian Zwiers\thanks{Email: \texttt{ian.zwiers@gmail.com}}
}
\date{}

\begin{document}

\maketitle

\begin{abstract}
In this note we study analytically and numerically the existence and stability of standing waves
for one dimensional nonlinear Schr\"odinger equations whose nonlinearities are the sum of three powers. Special attention is paid to the curves of non-existence and curves of stability change on the parameter planes.
\end{abstract}

\setcounter{tocdepth}{2}
{\hypersetup{linkcolor=brown}
\renewcommand{\baselinestretch}{0.9}\normalsize
\tableofcontents
\renewcommand{\baselinestretch}{1.0}\normalsize
}

\section{Introduction}
\label{sec1}

Consider the one dimensional nonlinear Schr\"odinger equations (NLS) for
$u(t,x):\R \times \R \to \CC$,
\EQ{
\label{NLS}
i\pd_t u +  \pd_x^2 u + f(u)=0,
}
with the nonlinearity $f(u):\CC \to \CC$ satisfying  $f(u)/u \to 0$ as $u \to 0$, and $f(\rho
e^{is} )e^{-is}= f(\rho)\in \R$ for any $\rho,s\in \R$.  A \emph{standing
wave} is a solution of \eqref{NLS} of the form  $u(t,x) = \phi(x) e^{i \om t}$ for some
$\om \in \R$ and a real-valued \emph{profile} $\phi \in L^2(\R)$, which then satisfies
\EQ{ \label{BL-eq}
\phi'' +f(\phi)=\om \phi.
}
We only consider solutions which decay rapidly at spatial infinity and hence assume $\om>0$.
The aim of this paper is
to examine the existence and stability of standing waves for nonlinearities
which are the sum of powers,
\EQ{
\label{f.def}
f(u) = \sum_{j=1}^m a_j |u|^{p_j-1} u, \quad 0\not=a_j \in \R, \quad 1<p_1<\cdots < p_m< \I.
}
The cases $m=1$ and $m=2$ (single and double power nonlinearities) are well studied and will be reviewed in Section \ref{sec2}. We will focus on $m=3$ in this paper with $p_j=j+1$.
Thus we
consider the triple power nonlinearity
\EQ{\label{234non}
f(u) =  a_1 |u|u + a_2 |u|^2 u+ a_3 |u|^3u,
}
with $a_1a_2a_3\not = 0$.
We may also consider other exponents and higher dimensions: A radial standing wave $u(x,t)=\phi(|x|) e^{i \om t}$ in $\R^n$ satisfies
\[
\phi''+ \frac{n-1}r \phi' + a_1|\phi|^{p_1-1}\phi + a_2 |\phi|^{p_2-1} \phi + a_3 |\phi|^{p_3-1}\phi=\om \phi.
\]
We limit ourselves to 1D with nonlinearity \eqref{234non} in this paper.

We now recall known results.
The existence and stability of standing wave solutions for general nonlinearities and general dimensions have a large literature. For existence, we refer to Berestycki and Lions \cite{BL83} and its references. Our Proposition \ref{th:BL83} is from \cite{BL83}. The standard references for stability are Grillakis, Shatah, and Strauss \cite{SS,GSS1,GSS2}. When restricted to one space dimension, one can decide the stability by the sign of a stability functional defined by an integral  found by Iliev and Kirchev  \cite{IK93}, see Proposition \ref{th1.2}.

We now consider standing waves for sums of power nonlinearities \eqref{f.def}.
For the single power nonlinearity $f(u)=|u|^{p-1}u$, the solutions are well known, and their stability problem has an extremely large literature and is still very active. Since it is not the concern of this paper, we only refer to the monograph 
\cite{Cazenave} by Cazenave for an introduction. 

Consider now a double power nonlinearity, $f(u)= a_1 |u|^{p_1-1}u + a_2 |u|^{p_2-1}u$. In the one dimensional setting, 
and if at least one of $a_j$ is positive, the standing waves $\phi_\om$ exist for $\om \in (0,\om^*)$ for
some $0<\om^* \le \infty$. If both $a_1$ and $a_2$ are negative, there is no standing wave.  
When $2p_1=p_2+1$, explicit formulas for the standing waves are known (see Appendix \S\ref{Explicit-formulas}). 
Their orbital stability is examined by Ohta \cite{Oh95}, and extended by Maeda \cite{Maeda}. See Proposition \ref{th1.3} for a summary.  
Note that the case $a_1<0<a_2$ with $p_2<5$ is not completely characterized yet. The small frequency case $\om \ll 1$ has been recently investigated by Fukaya and Hayashi \cite{Fukaya-Hayashi} for general dimensions.
For space dimension $n\ge2$, see Fukuizumi \cite{MR2029931},
Lewin and Nodari \cite{LN}, 
and Carles, Klein, and Sparber \cite{CKS}.
For related results,
see \cite{MR3449746,MR3927431,MR3976418} for the stability for cubic-quintic nonlinearity with an additional delta potential,
\cite{MR2387260} for the existence of standing waves for double power nonlinearity and harmonic potential,
and \cite{MR3556052,MR3409057} for blow up solutions for double power nonlinearities.

We are not aware of any study focused on the triple power nonlinearity, which is the subject of this paper.

Of great interest to us is the co-existence of stable and unstable standing waves for some fixed nonlinearities, with $\om$ being the only parameter. One wonders that, if a solution starts near a unstable standing wave, would it leave the neighborhood of unstable standing waves, and eventually converge to a \emph{stable} standing wave? This is partly motivated by our previous study \cite{NPT} and \cite{TVZ}. 
Another motivation is the study of the (in)stability of \emph{critical standing waves}. 
A standing wave family $\phi_{\om}e^{i\om t}$ may change from being stable to unstable as $\om$ goes across a critical $\om_c$.  It is shown by 
Comech and Pelinovsky \cite{MR1995870} that the critical $\phi_{\om_c}e^{i\om_c t}$ is unstable under certain conditions. Their results are extended by Ohta \cite{MR2785894} and
Maeda \cite{MR2923422} under various conditions. These results for NLS are extended to generalized KdV by Comech, Cuccagna, and Pelinovsky
\cite{CoCuPe}, and to derivative nonlinear Schr\"odinger equation by
Fukaya \cite{Fukaya},
Guo, Ning, and Wu \cite{MR4117079},
and Ning \cite{Ning}.
We hope that, for explicit nonlinearities in one space dimension, one is able to verify those conditions and do more detailed analysis.

We now consider our nonlinearity in more details.
 To study NLS \eqref{NLS} with $f(u)$ given by \eqref{234non},
we may let $u(x,t) = kv(\la^{-1} x, \la^{-2} t)$, $k,\la>0$. Then $v$ satisfies
\[
i\pd_t  v+ v_{xx} + \bar b |v|v +  \bar c |v|^2 v +  \bar d |v|^3v=0,
\]
\[
\bar b = a_1k\la^2, \quad 
\bar c = a_2k^2 \la^2, \quad
\bar d = a_3 k^3 \la^2.
\]
If we choose $k= |a_1/a_3|^{1/2}$ and $\la = |a_3/a_1^3|^{1/4}$, then $|\bar b|=|\bar d|=1$. Since $v$ and $u$ have the same qualitative properties,  we may assume $|a_1|=|a_3|=1$ without loss of generality. For the rest of this paper, we consider
\EQ{\label{234non2}
a_1 = \pm 1, \quad a_2 = -\ga, \quad a_3 = \pm 1,
}
and use $\ga \in \R$ as another parameter in addition to $\om$.  
To summarize, our standing wave profile $\phi:\R\to \R_+$ satisfies
\EQS{\label{phi234.eq}
\phi'' =g(\phi) = \om \phi - f(\phi),& \quad f(\phi)= a_1|\phi| \phi -\ga |\phi|^{2} \phi + a_3 |\phi|^{3}\phi,
\\[2mm]
\phi(t)>0 \quad \forall t\in\R,& \quad \lim _{t\to \pm \infty} \phi(t)=0,
}
with $a_1,a_3 = \pm 1$.
We have two parameters $\om>0$ and $\ga \in \R$, which lie on the half plane $(0,\I)\times \R$.

We consider 4 cases F*F, F*D, D*F and D*D, where the case F*D means $a_1=1$, $a_3=-1$, and consists of two subcases: FFD means $a_2>0$ and FDD means $a_2<0$.  The remaining cases are similarly defined. The negative sign in front of $\ga$ is convenient for the most interesting FDF case.  There are 8 combinations of the signs of the three terms. Except the case DDD, the remaining 7 cases have standing waves.

\begin{figure}[H]
\begin{minipage}[b]{.5\textwidth}
\centering
\includegraphics[width=1\textwidth, height=1\textwidth]{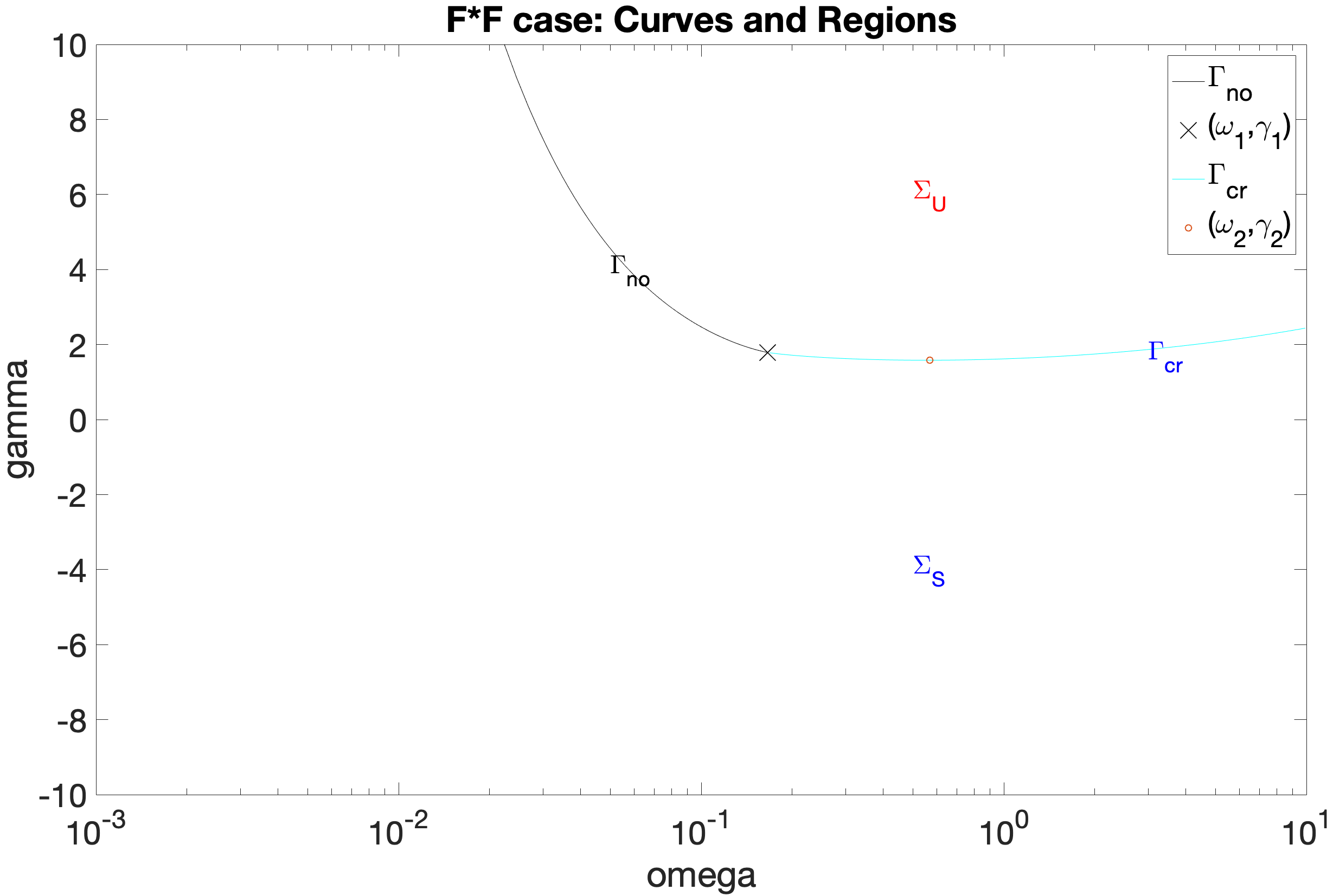}
\caption{(i) F*F case}
\label{log_region1}
\end{minipage}
\hfill
\begin{minipage}[b]{.5\textwidth}
\centering
\includegraphics[width=1\textwidth, height=1\textwidth]{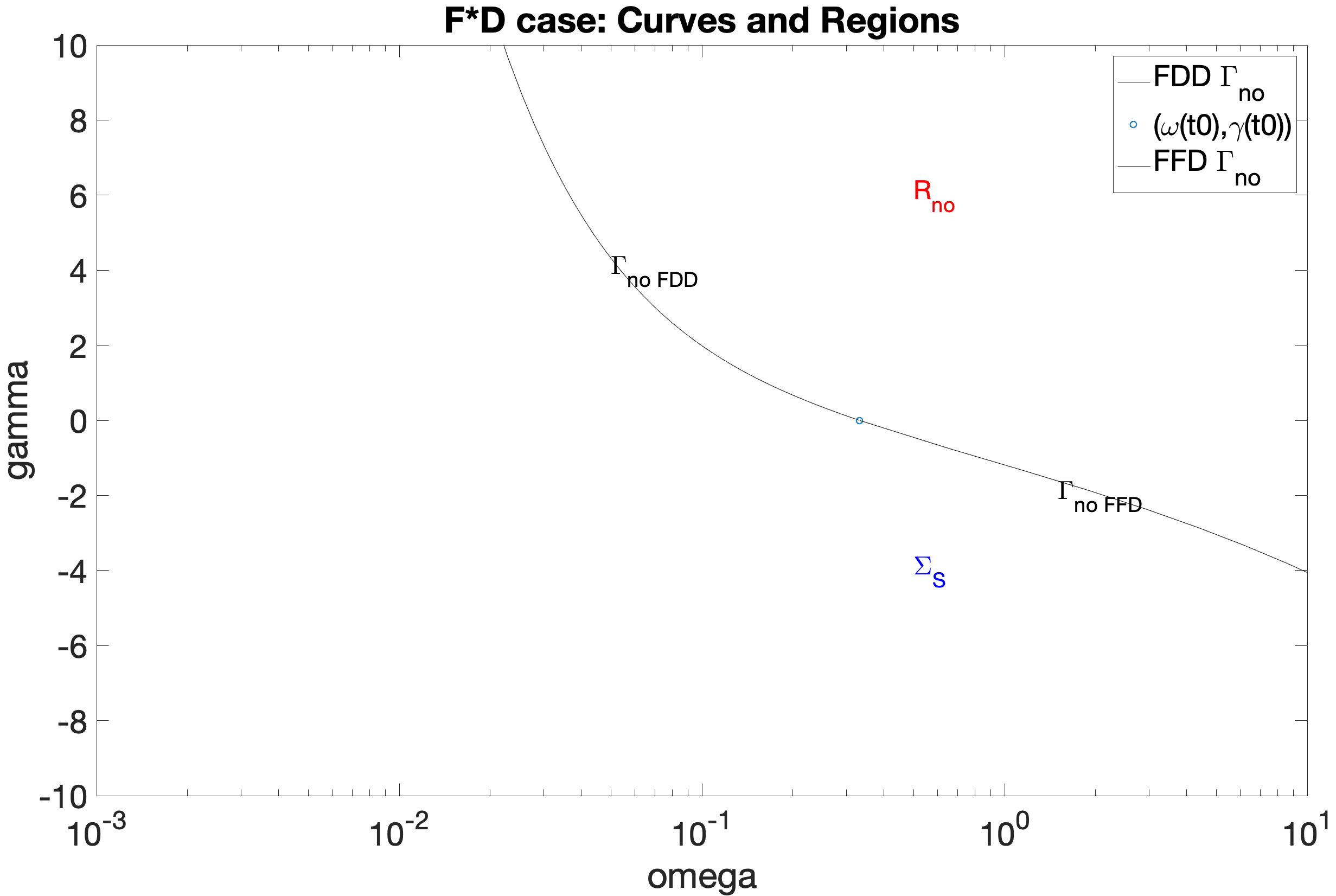} 
\caption{(ii) F*D case}
\label{log_region2}
\end{minipage}
\end{figure}

\begin{figure}[H]
\begin{minipage}[b]{.5\textwidth}
\centering
\includegraphics[width=1\textwidth, height=1\textwidth]{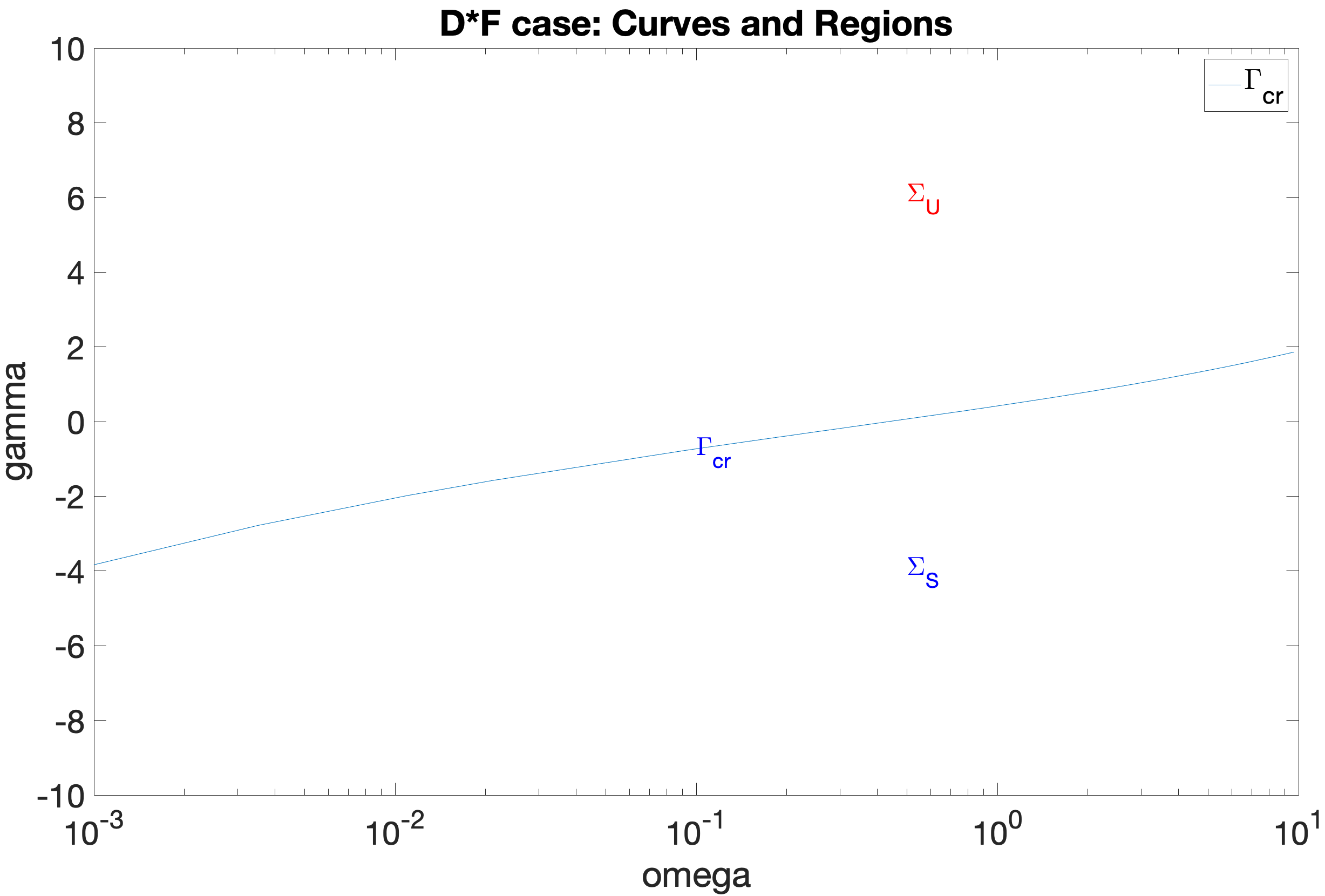}
\caption{(iii) D*F case}
\label{log_region3}
\end{minipage}
\hfill
\begin{minipage}[b]{.5\textwidth}
\centering
\includegraphics[width=1\textwidth, height=1\textwidth]{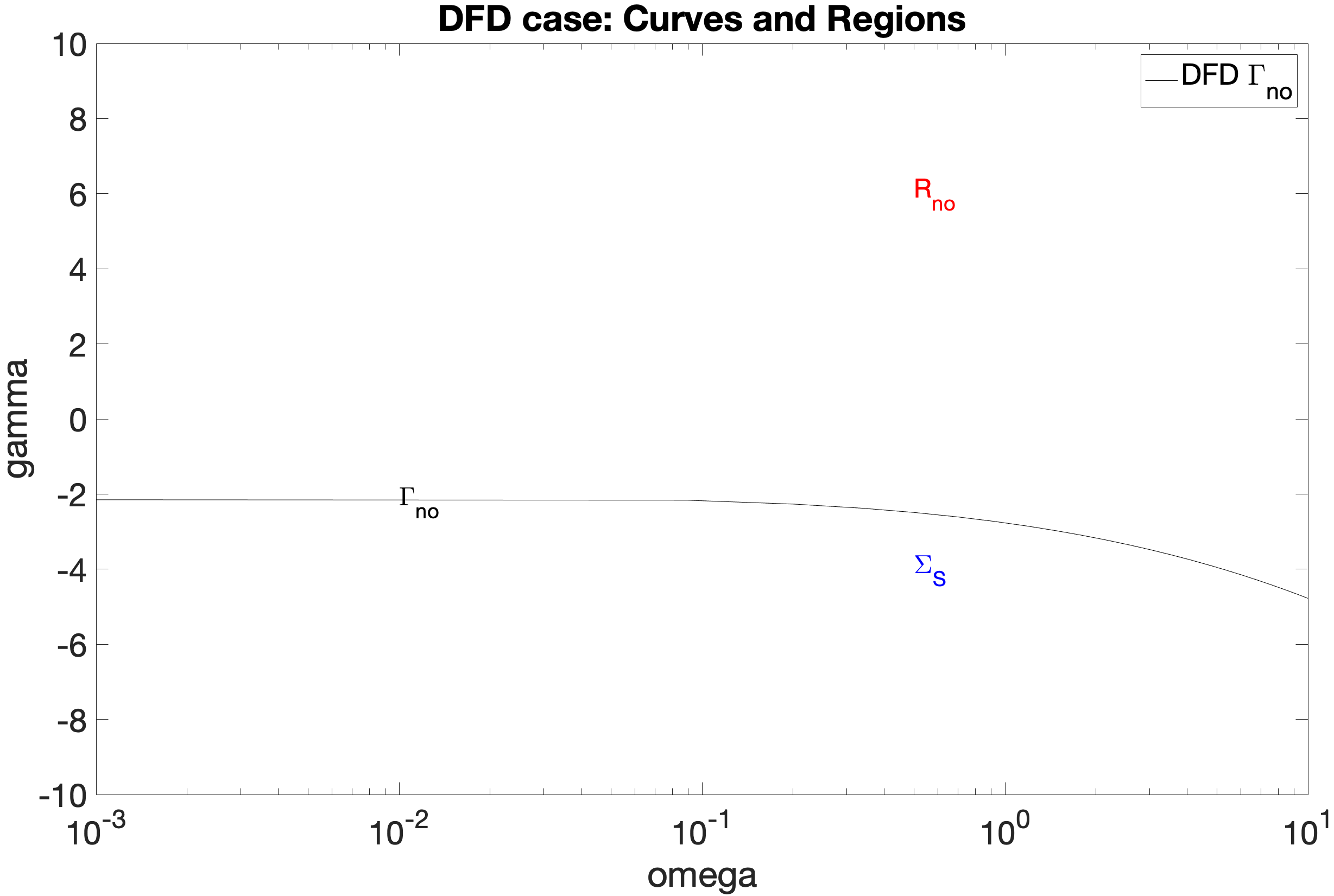} 
\caption{(iv) D*D case}
\label{log_region4}
\end{minipage}
\end{figure}

We now present our numerically computed non-existence and stability regions on the parameter half-planes for the 4 cases. See Figures \ref{log_region1}--\ref{log_region4}. The parameter $\om$ is presented in log scale to zoom in the small $\om$ part.

In the figures,
\begin{enumerate}
\item $R_{\no}$ is the set of $(\om,\ga)$ that $\phi$ does not exist, present in (ii) and (iv);
\item $\Si_{S}$ is the set of $(\om,\ga)$ that $\phi$ exists and is stable, present in all cases;
\item $\Si_{U}$ is the set of $(\om,\ga)$ that $\phi$ exists and is unstable, present in (i) and (iii);
\item $\Ga_{\no}$ is the curve of $(\om,\ga)$ that $\phi$ does not exist, present in (i), (ii) and (iv);
\item $\Ga_{cr}$ is the curve of $(\om,\ga)$ that $\phi$ exists and is at threshold of stability, present in (i) and (iii).
\end{enumerate}

The boundary curves $\Ga_{\no}$ for existence will be analytically computed in Section \ref{sec3}, based on the study of double zeros of the potential function $G(x)=\int_0^x g(s)ds$ for $g(s)$ given in \eqref{phi234.eq}.
In contrast, the curves $\Ga_{cr}$ for stability change are computed only numerically, see Section \ref{sec4}. 

Among the 4 cases, the F*F case seems the most interesting, as $\Ga_{\no}$ and $\Ga_{cr}$ co-exist and meet at
\[
(\om_1,\ga_1):=(\frac{2\sqrt 5}{27}, \frac 4{\sqrt 5}) \approx (0.1656,1.7889).
\]
Moreover, $\Ga_{cr}$ has a minimal $\ga$ value at (see \eqref{best-approx})
\[
(\om_2,\ga_2) \approx (\om_2^{**},\ga_2^{**}) =(0.554837092755109,1.58170475989899).
\]
The pair $(\om_1,\ga_1)$ is computed analytically, while $(\om_2,\ga_2)$ is only computed numerically.
The accuracy of $\gamma_2^{**}$ is about $10^{-9}$, but the accuracy of $\om_2^{**}$ is about $10^{-4}$, much larger. See \S\ref{sec5.2}.

We now explain the structure of the paper. In Section \ref{sec2}, we recall the general existence and stability results, and the known results for double power nonlinearities.

 In Section \ref{sec3}, we analyze the existence for the triple power nonlinearity \eqref{234non}-\eqref{234non2} by studying the double zeros of the potential function $G(x;\om,\ga)$. We characterize the existence regions and their boundary $\Ga_\no$. We also study solutions inside the standing wave homoclinic orbit in Subsection \ref{sec3.4}.
 
In Section \ref{sec4}, we study the stability regions by studying the stability function $J(\om,\ga)$. We study analytically the limits of $J(\om,\ga)$ as $(\om,\ga)$ approaches the non-existence curve $\Ga_{\no}$ from different sides in  Subsection \ref{sec4.1}, and characterize the  stability regions for the four cases using the level sets of $J$ in  Subsection \ref{sec4.2}.
  
In Section \ref{sec5} we discuss our numerical methods and observations.
In Subsection \ref{sec5.1}, we describe the computations of the level curves of the stability functional $J(\omega,\gamma)$ and the stability regions for the four cases.
In Subsection \ref{sec5.2}, we describe the computation of the curve of stability change $\Gamma_{cr}$ and its minimal point $ (\omega_2, \gamma_2)$ in the F*F case.
In Subsection \ref{sec5.4}, we describe the 3 methods that we use for computing the standing wave $\phi_{\omega,\gamma}$.

In Appendix Section \ref{Explicit-formulas}, we give alternative explicit formulas for the standing waves of \eqref{35phi.eq} for double power nonlinearities with $2p_1 = p_2+1$.

\section{Preliminaries}
\label{sec2}

In this section we recall some general results for bounded solutions $u:\R \to \R$ of %
\EQ{\label{g=omu-f}
u''=g(u)= \om u - f(u), \quad \om >0, \quad \lim_{u \to 0} \frac {f(u)}u=0.
}

\subsection{Existence}
\label{sec2.1}

The following is a general existence result.

\begin{proposition}[Existence \cite{BL83}]\label{th:BL83}
Let $g\in C(\R;\R)$ be a locally Lipschitz continuous function with
$g(0)=0$ and let $G(t)=\int_0^t g(s)\,ds$. A necessary and sufficient condition for the existence of a solution $\phi$ of the problem
\begin{align}
\phi \in C^2(\R), \quad    \lim_{t \to \pm \I}\phi(t)=0,& \quad
\phi(0)>0,
\\
\phi''=g(\phi),& \label{phi.eq}
\end{align}
is that
\begin{equation} 
\label{eq1.4}
\phi_0 =\inf \{ t>0:\ G(t)=0\} \text{ exists},\quad \phi_0>0,
\quad g(\phi_0)<0.
\end{equation}
\end{proposition}

For the special case \eqref{g=omu-f}, the existence can be derived directly from the phase plane analysis for \eqref{phi.eq} without use of Proposition \ref{th:BL83}. We now describe it. The corresponding planar first order system for \eqref{phi.eq} is, with $x=\phi$ and $y=\dot \phi$,
\EQ{\label{planar}
\dot x = y, \quad \dot y = g(x).
}
Every solution moves on a level curve of the total energy
\[
E(x,y) = \frac 12 y^2 - G(x) =   \frac 12 y^2 -   \frac \om 2 x^2 + F(x),
\]
where $F(x) = \int_0^x f(t)\,dt$.
Assume that $G$ has a smallest positive zero $\phi_0>0$, $G(\phi_0)=0$. We have $g(\phi_0) \le 0$ because $G(x)>0$ for $0<x<\phi_0$.
As $E(x,0)=-G(x)<0$ for small $x>0$, by continuity, 
\[
E(x,0) =-G(x)<0, \quad 0<x<\phi_0.
\]
For the single power nonlinearity $f(u) = |u|^{p-1}u$ and $g(u) = \om u - f(u)$, the solution passing through $(x,y)=(\al,0)$ with 
 $0<\alpha <\phi_0$ is either a fixed point or a periodic orbit. But this is not true for general $f(u)$. 
See Lemma \ref{th3.3} and Examples \ref{ex3.4}--\ref{ex3.5} for the triple power case.

\medskip

The nature of the solution passing through $(x,y)=(\phi_0,0)$ depends on whether $g(\phi_0)=0$. Let $\Om$ denote the connected component of the sublevel set $E<0$ that contains the open line segment from $(0,0)$ to $(\phi_0,0)$. It is bounded, and $E=0$ on its boundary $\pd \Om$. We assume $\pd\Om$ is a nice curve lying in between  $0\le x \le \phi_0$. Its fixed point(s) lie on the $x$-axis since $\dot x =y\not =0$ if $y\not=0$.
\EN{
\item If $g(\phi_0)=0$: $(\phi_0,0)$ is a fixed point of \eqref{planar}. The upper branch of $\pd\Om$ in the first quadrant is 
 a heteroclinic orbit from $(0,0)$ to $(\phi_0,0)$.
\item If $g(\phi_0)<0$: $(\phi_0,0)$ is not a fixed point. The curve $\pd\Om$ is a homoclinic orbit from $(0,0)$ to $(0,0)$ passing through $(\phi_0,0)$. The $x$-component of this orbit is always positive.
}

Formally the solution passing $(\phi_0,0)$ satisfies $E= \frac 12 y^2 - G(x)=0$. 
Suppose $(\phi_0,0)=(x(0),y(0))$. 
Then the branch in the fourth quadrant (corresponding to $0<t<\infty$) satisfies
\EQ{\label{dt.formula}
\frac {dx}{dt} = -\sqrt{2G(x)}, \quad dt = -\frac {dx}{ \sqrt{2G(x)}}.
}
If $x(0)=\phi_0$ and $x(t)=\phi$, then
\EQ{\label{t.formula}
t = -\int_{\phi_0}^\phi \frac {dx}{ \sqrt{2G(x)}}.
}
If $g(\phi_0)=0$, then $\phi_0$ is a double zero of $G$, and the integral \eqref{t.formula} is not integrable. This agrees with the fact that $(\phi_0,0)$ is a fixed point. If $g(\phi_0)<0$, then \eqref{t.formula} is a well-defined improper integral. The integrand is positive for any $x \in (0,\phi_0)$, and hence the integral is well-defined. As $x \to 0_+$, $\sqrt{2G(x)} \sim \sqrt \om x$. Hence $t \to \infty$ as $\phi \to 0_+$. The inverse function of the integral function \eqref{t.formula} is our desired solution $\phi(t)=\phi_\om(t)$.

We may now treat $\om$ as a parameter and denote $\phi=\phi_\om$.
The set of $\om$ for which $\phi_\om$ exists is open because the existence condition \eqref{eq1.4} that $\phi_0$ is that first zero of $G$ and that $g(\phi_0)<0$ are
preserved under small perturbations in $\om$. %

\subsection{Stability}
\label{sec2.2}

The \emph{orbit} of a standing wave $\phi e^{i \om t}$  is the
set  obtained from it by translation and phase shift,
\[
\cO_\phi= \{  e^{i\eta}\phi(\cdot-\xi): \quad (\xi,\eta)\in \R^2\} .
\]
The distance of $u$ to this orbit of $\phi$ is
\[
\text{dist}(u,\cO_\phi)=\inf_{(\xi,\eta)\in \R^2} \norm{u - e^{i\eta}\phi(\cdot-\xi)}_{H^1(\R)}.
\]
We say a standing wave $\phi e^{i \om t}$ is \emph{orbitally stable} if
\[
\forall \e>0,\ \exists \de>0, s.t.\ \text{dist}(u(0),\cO_\phi) < \de \Rightarrow
\text{dist}(u(t),\cO_\phi)<\e \ \forall t \ge 0,
\]
where $u(t)$ is the solution of \eqref{NLS} with initial data $u(0)$.
This definition
falls in the general framework of \cite{GSS2}, and differs from \cite{SS,GSS1} since the orbit contains translations. We may remove the translations from $\cO_\phi$ if we restrict our perturbations to even perturbations so that the solutions do not move. For general (non-even) perturbations, 
the solutions may get boost from the perturbations and start to move, and hence
we need to contain translations in $\cO_\phi$. 
For example, we may get a traveling wave by the 
Galilean transformation for $v\not =0$,
\[
u(x,t)= \phi_\om(x-vt) e^{i (\om t + vx/2 - tv^2/4)}. 
\]

We prepare a few definitions before we state an orbital stability result.
For a family of standing waves $\phi_\om$, $\om \in (\om_a,\om_b)$,
denote
\[
d(\om)=E(\phi_\om)+\om Q(\phi_\om),
\]
where 
\[
E(u)=\int_\R |u_x|^2 - 2F(u)\,dx, \quad Q(u) = \int_\R |u|^2 dx.
\]
Then, assuming enough regularity of $f$, $d(\om)$ is $C^2$ and
\[
d'(\om) = Q(\phi_\om).
\]

Following \cite{IK93}, we denote for $s \ge 0$
\EQS{\label{IK.def}
f(s)&=- f_1(s^2)s, \quad F_1(s) = \int_0^s f_1(t)dt=-2F(\sqrt s),
\\
U(s) &= \om s + F_1(s) = \om s - 2F(\sqrt s) = 2G(\sqrt s), \quad
U'(s)= \frac{g(\sqrt s)}{\sqrt s}.
}
($f_1$ and $F_1$ are denoted as $f$ and $g$ in  \cite{IK93}.)
The existence condition \eqref{eq1.4} is equivalent to (with $a=\phi_0^2$)
\begin{equation} 
\label{eq1.6}
\exists a\in (0,\I)  \text{ such that } U(a)=0,
\quad U'(a)<0, \quad U(s)>0\quad (0<s<a).
\end{equation}

The following is an orbital stability result by Iliev and Kirchev  \cite{IK93}.

\begin{proposition}[Orbital stability \cite{IK93}] 
\label{th1.2}
Suppose $f(u)$ is such that \eqref{NLS} is locally wellposed in $H^2(\R)$,  there is a constant $A>0$ such that $f_1(s) \in C^0[0,A) \cap
  C^1(0,A)$, $sf_1'(s) \to 0$ as $s \to 0$, and \eqref{eq1.6} is satisfied
  with $a<A$. If $d''(\om)>0$,
  then $\phi_\om e^{i\om t}$ is orbitally stable. If $d''(\om)<0$, then
  $\phi_\om e^{i\om t}$ is orbitally unstable.
Furthermore, 
\EQ{\label{th1.2-1}
d''(\om)=\frac d{d\om}\int \phi_\om^2(x)dx =  \frac {-1}{2U'(a)}\int_0^a \bke{3+\frac
      {as[f_1(a)-f_1(s)]}{aF_1(s)-sF_1(a)}}\bke{\frac s{U(s)}}^{1/2}  ds.
}
\end{proposition}

The last formula is \cite[Lemma 6]{IK93}.
Note that  $ f_1(a) - f_1(s)= U'(a) - U'(s)$ and
\[
aF_1(s)-sF_1(a) = a(U(s)-\om s) - s (U(a)-\om a) = aU(s).
\]
Hence
\EQ{\label{th1.2-1b}
d''(\om)=  \frac {-1}{2U'(a)} \int_0^a \bke{3 + \frac {s(U'(a)-U'(s))}{U(s)}}
\bke{\frac s{U(s)}}^{1/2}  ds.
}
Changing variable $s=a\si$ and relabeling $\si$ as $s$, we get the form we use in MATLAB:
\EQ{
\label{eq1.20.1}
d''(\om)=
\frac {-a^{1.5}}{2U'(a)}\int_0^1
\bke{3 + \frac {as(U'(a)-U'(as))} {U(as)}}\frac {\sqrt{s}} {\sqrt{U(as)}}\,ds.
}
This formula is convenient for numerics since it has a fixed domain $(0,1)$ and we only need to compute the vectors $U_k=U(ak/N)$ and $V_k=U'(ak/N)$, $0\le k \le N$, and do not need to compute $f_1$ or $F_1$.

If we make the natural change of variables $s=x^2$ (with $x=\phi_\om(t)$), using \eqref{IK.def} and $a=\phi_0^2$), we get
\EQ{\label{th1.2-1c}
d''(\om) = \frac{-\sqrt 2}{4g(\phi_0)} \int_0^{\phi_0} \bket{ \frac {6 \phi_0}{ \sqrt{G(x)}}  + \frac {x[x g(\phi_0) -\phi_0 g(x)  ]}{ \sqrt{G(x)}^3}  }\,x^2dx .
}
Note that we are free to add into the integral of \eqref{th1.2-1c}
\EQ{
0=  \int_0^{\phi_0} \frac{d}{dx}\bke{ \frac {2A(x)}{ \sqrt{G(x)}}} dx = 
 \int_0^{\phi_0}\bke{ \frac {2A'(x)}{ \sqrt{G(x)}}  -  \frac {A(x) g(x)}{ \sqrt{G(x)}^3}} dx 
}
for any $A(x) \in C^1([0,\phi_0])$ with $A(0)=A'(0)=A(\phi_0)=0$. For example, the choice $A(x)=-\phi_0x^3$ does not satisfy $A(\phi_0)=0$, while  the choices $A(x)=x^4-\phi_0x^3$ and $A(x)=\phi_0^2x^2-\phi_0x^3$ are valid and give us
\EQ{\label{th1.2-1d}
d''(\om) = \frac{-\sqrt 2}{4g(\phi_0)} \int_0^{\phi_0} \frac { 8x^3 G(x)+ x^4 g(\phi_0) -x^4 g(x) }{ \sqrt{G(x)}^3}  \,dx
}
and
\EQ{
d''(\om) = \frac{-\sqrt 2}{4g(\phi_0)} \int_0^{\phi_0} \frac { 4\phi_0^2x G(x)+ x^4 g(\phi_0) -\phi_0^2x^2 g(x) }{ \sqrt{G(x)}^3}  \,dx.
}

\subsection{Sum of positive powers}
\label{Sum-of-positive-powers}
We now consider the special case that $f(u)$ is a sum of $m$ positive power nonlinearities, given as in  \eqref{f.def}, with all coefficients $a_i>0$.

In the classical case $m=1$ and $f(u) = a_1 |u|^{p_1-1}u$, a lot is known: We need $a_1>0$ to ensure the existence.
In this case $\phi_\om$ exist for all $\om$, and indeed are the rescaling of each other. 
It is well known that they are stable if $p_1 < 5 $ and unstable if $5 \le p_1$. This result follows from Proposition \ref{th1.2} if $p_1 \not =5$ (see the argument for $m>1$ case below), and needs extra work if $p_1=5$.  See e.g.~\cite{Cazenave}.

When $m>1$, $f(u)$ is given by \eqref{f.def}, and all $a_1,\cdots, a_m$ are positive, we have
\[
G(u)= \frac \om 2 |u|^2 - F(u), \quad F(u)=\sum_{k=1}^m \frac {a_k}{p+1} |u|^{p+1}.
\]
It clearly has a first positive zero $\phi_0$. 
Moreover,
\[
\frac12\phi_0 g(\phi_0) = \frac \om 2 \phi_0^2 - \sum_{k=1}^m \frac {a_k}{2} \phi_0^{p+1}<G(\phi_0)=0.
\]
Thus $g(\phi_0)<0$.
By Proposition \ref{th:BL83}, the family of standing waves $\phi_\om$ exist for all $0<\om<\I$. 

For their stability, rewrite \eqref{th1.2-1b} as 
\EQ{
\label{eq1.20}
d''(\om)=
\frac {-1}{2U'(a)}\int_0^a
I(s)\bke{\frac s{U(s)}}^{\frac32} ds,
\quad
I(s) =  \frac{3U(s)}s + f_1(a) - f_1(s) .
}
Using $U(s)/s=2G(\sqrt s)/s =\om - \sum_{k=1}^m \frac {2a_k}{p_k+1} s^{\sigma_k}$, where $\sigma_k=\frac{p_k-1}2$, $U(a)=0$ hence $\om=\sum_{k=1}^m \frac {2a_k}{p_k+1} a^{\sigma_k}$,  and $f_1(s)=-f(\sqrt s)/\sqrt s = - \sum_{k=1}^m a_k s^{\frac{p_k-1}2}$, we get
\EQN{
I(s)&=\frac{3U(s)}s + f_1(a) - f_1(s)
\\
&
= \sum_k \frac{6a_k}{p_k+1} (a^{\sigma_k}-s^{\sigma_k})-\sum_k a_k (a^{\sigma_k}-s^{\sigma_k})
 \\
&= \sum_k \frac{a_k(5-p_k)}{p_k+1} (a^{\sigma_k}-s^{\sigma_k}).
}
Since $m>1$, $I(s)>0$ for $0<s<a$ if $p_m \le 5$. In this case, $\phi_\om$ for all $\om$ are orbital stable by Proposition \ref{th1.2} and \eqref{eq1.20}. On the other hand, if $p_1  \ge 5$, $I(s)<0$ for $0<s<a$, hence $\phi_\om$ for all $\om$ are orbital unstable. The intermediate cases $p_1<5 <p_m$ are more subtle.

\subsection{Double power nonlinearities}
\label{sec2.4}

When $m=2$ and $f$ is a double power nonlinearity,
\EQ{
\label{2power}
f(u)= a_1 |u|^{p_1-1}u + a_2 |u|^{p_2-1}u, \quad 1<p_1 <p_2 < \I, \quad a_1a_2\not =0,
}
and if at least one of $a_j$ is positive,
the  standing waves $\phi_\om$ exist for $\om \in (0,\om^*)$ for
some $0<\om^* \le \infty$. If both $a_1$ and $a_2$ are negative, there is no standing wave.  Their stability 
is examined by Ohta \cite{Oh95}, and extended by Maeda \cite{Maeda}. To present their results in a compact form, we make the following

\medskip

{\bf Definition.}\
The family $\bket{\phi_\om}_{0<\om<\om^*}$ is of type SU if there exist $0<\om_1<\om^*$ so that
$\phi_\om$ is stable for $\om \in (0,\om_1)$, 
and unstable 
for  $\om \in (\om_1,\om^*)$. 
It is of type U?S if there exist $0<\om_1<\om_2<\om^*$ so that
$\phi_\om$ is unstable for $\om \in (0,\om_1)$, unknown for  $\om \in (\om_1,\om_2)$,
and stable 
for  $\om \in (\om_2,\om^*)$. There is no assertion at $\om_1$ and $\om_2$.
Other types are defined similarly.

\medskip

We summarize the results of \cite{Oh95,Maeda} in the following.

\begin{proposition}
\label{th1.3}
Let $f(u)$ be of the form \eqref{2power}. We consider the stability type of the family
$\bket{\phi_\om}_{0<\om<\om^*}$ in three groups: %

\medskip

\noindent
(focusing-focusing, FF) Let $a_1,a_2>0$. Then $\om^*=\I$.

(a) If $p_2 \le 5$, it is of type S.

(b) If $p_1 \ge 5$, it is of type U.

(c) If $p_1<5<p_2$, it is of type SU.

\medskip

\noindent(focusing-defocusing, FD) Let $a_1>0$, $a_2<0$. Then $\om^*<\I$.

(a) If $p_1 \le 5$, it is of type S.

(b) If $p_1 > 5$, it is of type US.

\medskip
\noindent(defocusing-focusing, DF) Let $a_1<0$, $a_2>0$. Then $\om^*=\I$.

(a) If $p_2 \ge 5$, it is of type U.

(b) If $p_2<5$, it is of type ?S.

\quad\ \
If $p_2<5$ and $p_1+p_2>6$, it is of type U?S.

\quad\ \
If $p_2<5$, $p_1+p_2>6$, and $7/3<p_1$, it is of type US.

\quad\ \
If $p_1=2$ and $p_2=3$, it is of type S.

\end{proposition}

Note that the case (DF)-(b) is not completely decided. However, it is conjectured that
the stability change occurs at most once \cite[p.265]{Maeda}. The instability of standing waves with small frequency $\om$ in the DF case
has been recently investigated by Fukaya and Hayashi \cite{Fukaya-Hayashi} for general dimensions. For dimension 1, standing waves with sufficiently small $\om$ are unstable if $1<p_1<p_2<5$ and $(p_1+3)(p_2+3)>32$.

Note that, for double power nonlinearity with $2p_1=p_2+1$, explicit formulas for the standing waves are known. We will give a new formulation which gives the standing waves for all FF, FD and DF cases in one single formula in Appendix \S\ref{Explicit-formulas}.

\section{Existence for triple power nonlinearities}
\label{sec3}

In this section, we study the existence of the standing wave profile $\phi:\R\to \R_+$ satisfying \eqref{phi234.eq} with triple power nonlinearity for parameters $\om>0$ and $\ga\in\R$, that is,
\EQN{
\phi'' =g(\phi) = \om \phi - f(\phi),& \quad f(\phi)= a_1|\phi| \phi -\ga |\phi|^{2} \phi + a_3 |\phi|^{3}\phi,
\\
\phi(t)>0 \quad \forall t\in\R,& \quad \lim _{t\to \pm \infty} \phi(t)=0,
}
for the four cases $a_1,a_3 = \pm 1$. Recall %
$G(x) = \int_0^x g$. %
For our special choice of $f$,%
\EQ{\label{G234.def}
G(x) = \frac \om 2 x^2 -  \frac{a_1}3|x|^3 +\frac{\ga}4 x^{4}  -\frac{ a_3}5 |x|^{5}, \quad (x \in \R).
}

\subsection{Definition and basic properties of existence regions of parameters}

By Section \ref{sec2.1}, the standing wave profile $\phi$ exists if and only 
 if that 
\begin{enumerate}
\item [(E1)]  the first positive zero $\phi_0$ of $G(x)$ exists, and
\item [(E2)] $g(\phi_0)<0$.
\end{enumerate}
(E1) fails if $G(x)>0$ for all $x>0$. When (E1) is valid, (E2) fails if $\phi_0$ is a double zero of $G$, $g(\phi_0)=0$.

Condition (E2) is preserved under small perturbations of $\om$ and $\ga$. Condition (E1) by itself is not so, because the first positive zero  may disappear or jump under perturbations of $\om$ and $\ga$. However, when any of these happens, it is necessary that $g(\phi_0)=0$, violating (E2). We formulate this as  a lemma.

\begin{lemma} \label{th3.1}
Suppose (E1) and (E2) are valid at $(\om_0,\ga_0)\in \R_+\times \R$. Then it is valid in a small neighborhood of $(\om_0,\ga_0)$. Moreover, the map $(\om,\ga)\mapsto \phi_0(\om,\ga)$ is  continuous in the neighborhood.
\end{lemma}
\begin{proof} Let $\phi_0$ be the first zero of $G$ at $(\om_0,\ga_0)$. Let
$H(x,\om,\ga)=\frac \om 2 x^2 -  \frac{a_1}3x^3 +\frac{\ga}4 x^{4}  -\frac{ a_3}5 x^{5}$. We have
$H(\phi_0,\om_0,\ga_0)=0$ and $H_x(\phi_0,\om_0,\ga_0)<0$. By the implicit function theorem, there is a continuous map $X(\om,\ga)$ in a small neighborhood of $(\om_0,\ga_0)$ such that $H(X(\om,\ga), \om,\ga)=0$ and $X(\om_0,\ga_0)=\phi_0$. We also have $H_x(X(\om,\ga), \om,\ga)<0$ by continuity. We may choose a smaller neighborhood to ensure that $X(\om,\ga)$ remains the first positive zero.
\end{proof}

For each of the 4 cases $a_1,a_3 = \pm 1$, %
we denote by $R_{\ex}$ the subset of parameters $(\om,\ga)\in \R_+\times \R$ for which the standing wave profile $\phi$ exists. That is, both (E1) and (E2) are valid. By Lemma \ref{th3.1}, $R_{\ex}$ is an open subset of $\R_+\times \R$. Let $R_{\no}$ be the complement of $R_{\ex}$. It is a closed subset of $\R_+\times \R$. Let $\Ga_\no$ be the boundary of $R_\no$ in $\R_+\times \R$.

\begin{lemma} \label{th3.1b}
Consider the first positive zero $\phi_0$ of $G$ as a function of $(\om ,\ga)\in R_\ex$. It is continuously differentiable and
$\frac {\pd \phi_0}{\pd \om}>0$ and $ \frac {\pd \phi_0}{\pd \ga}>0$ for all  $(\om ,\ga)\in R_\ex$ and for all 4 cases.
\end{lemma}

\begin{proof} Let $H(x,\om,\ga)=\frac \om 2 x^2 -  \frac{a_1}3x^3 +\frac{\ga}4 x^{4}  -\frac{ a_3}5 x^{5}$. We have $H(\phi_0(\om,\ga),\om,\ga)=0$ and $g(\phi_0(\om,\ga))<0$. Hence
\[
0 = \pd_\om H(\phi_0(\om,\ga),\om,\ga)=\frac 1 2 \phi_0^2 + g(\phi_0) \frac {\pd \phi_0}{\pd \om},\quad 
0 = \pd_\ga H(\phi_0(\om,\ga),\om,\ga)=\frac 1 4 \phi_0^4 + g(\phi_0) \frac {\pd \phi_0}{\pd \ga}. 
\]
Thus $\frac {\pd \phi_0}{\pd \om} = \frac { \phi_0^2}{-2g(\phi_0)}>0$ and $\frac {\pd \phi_0}{\pd \ga} = \frac { \phi_0^4}{-4g(\phi_0)}>0$.
\end{proof}

\begin{lemma} \label{th3.2}
At any $(\om_0,\ga_0)\in \Ga_\no$, (E1) is valid and (E2) fails. That is, the first positive zero of $G$ at $(\om_0,\ga_0)$ exists and is a double zero. 
\end{lemma}
\begin{proof} Choose $(\om_k, \ga_k) \in R_\ex$, $k\in \N$, with
$(\om_k, \ga_k)  \to (\om_0,\ga_0)$ as $k \to \infty$. Let $\phi_0^k=\phi_0(\om_k, \ga_k)$, their first positive zero of $G$.
Since
\[
\frac {\om_k} 2  -  \frac{a_1}3x +\frac{\ga_k}4 x^{2}  -\frac{ a_3}5 x^{3}=0,\quad x=\phi_0^k,
\]
$\phi_0^k$ is uniformed bounded away from zero and infinity. There is a subsequence, still denoted as  $\phi_0^k$, that converges to some finite $x_0>0$. Taking limits of $H(\phi_0^k,\om_k, \ga_k) =0$, we get $H(x_0,\om_0,\ga_0)=0$.
Thus $G$ has a positive zero $x_0$ at $(\om_0,\ga_0)$ and (E1) is valid. By assumption, (E2) fails. Either $x_0$ is a double zero, or the first positive zero $\phi_0$ is less than $x_0$. In the latter case, $\phi_0$, $H(\phi_0,\om_k, \ga_k)>0$ and $\lim_{k \to \infty} H(\phi_0,\om_k, \ga_k)=0$. Hence $\phi_0$ is a double zero.
\end{proof}

\subsection{Double zeros of the potential function $G$}
\label{S:00}
Suppose (E1) is valid so that $G$ has a positive zero. Then (E2) is valid if and only if the first positive zero of $G$ is not a double zero.
Let us now consider all positive double zeros $t$ of $G$ (not necessarily the first positive zero),
\[
G(t)=0, \quad G'(t)=g(t)=0, \quad t>0.
\]
They correspond to fixed points $(t,0)$ of \eqref{planar} with zero energy, $E(t,0)=0$. Since $G(x)$ is a polynomial of degree $5$ and both $x=0$ and $x=t$ are its double zeros,
\EQ{\label{G.factorization}
G(x) = - \frac{ a_3}5 x^2(x-t)^2 (x-x_0)
}
for some $x_0 \in \R$, $x_0\not=0$. Thus for every $(\om,\ga)$ there is at most one positive double zero $t$. When there is one,
solving $2G(t)/t^2 = g(t)/t=0$ using \eqref{G234.def}, we get
\EQ{\label{om-c.formula}
\om = a_1\frac t3 - a_3\frac 15 t^3, \quad \ga=a_1\frac 2{3t} +  a_3\frac 65 t ,
\quad x_0 =  \frac{5 \om}{2a_3\, t^2} = \frac{5a_1}{6a_3\, t} - \frac 12 t  .
}
This gives a candidate curve of parameters for nonexistence. 
The requirement $\om>0$ gives
\EQ{\label{om>0}
\frac 53 a_1 > a_3t^2, \quad a_3 x_0>0.
}
If $x_0 \not \in (0, t)$,
the double zero $t$ is the first positive zero, and $\phi_\om$ does not exist.
When $x_0 \in (0,t)$, then $x_0$ is the first positive zero and $G'(x_0)<0$. Thus $\phi_\om$ exists with $\phi_\om(0)=x_0$. 

Note that, whenever $d\om/dt \not =0$,
\EQ{\label{negative-slope}
\frac{d\ga}{d\om} = \frac{d\ga/dt}{d\om/dt}  = -\frac2{t^2}<0.
}

\subsection{Characterization of non-existence regions of parameters}

We now identify the non-existence curve $\Ga_\no$ and non-existence region $R_\no$ for each of the four cases $a_1,a_3=\pm1$.

\begin{theorem}[Characterization of non-existence regions]\label{th3.4}

The non-existence curve $\Ga_\no$ and non-existence region $R_\no$ 
of parameters $(\om,\ga)\in\R_+\times\R$ for each of the four cases are:
\begin{enumerate}
\item F*F case $a_1=a_3=1$:  The curve $R_{\no}=\Ga _{\no}$ is given by 
\EQ{\label{GammanoFF}
\om = \frac t3 - \frac 15 t^3, \quad \ga=\frac 2{3t} + \frac 65 t,
\quad 0<t \le  \frac{\sqrt{5}}3.
}
We have $(\om,\ga) \to (0,\infty)$ as $t \to 0_+$ and $(\om,\ga) =(\om_1,\ga_1) =(\frac {2\sqrt5}{27},\frac{4\sqrt5}5) $ as $t =\frac{\sqrt{5}}3$.

\item F*D case $a_1=1=-a_3$:  $R_{\no}$ is the set on and above the curve $\Ga _{\no}$ given by
\EQ{\label{GammanoFD}
\om = \frac t3 + \frac 15 t^3, \quad \ga=\frac 2{3t} - \frac 65 t,
\quad 0<t <  \infty.
}
We have $(\om,\ga) \to (0,\infty)$ as $t \to 0_+$ and $(\om,\ga)\to (\infty,-\infty)$ as $t \to  \infty $.

\item D*F case $a_3=1=-a_1$: Both $R_{\no}$ and $\Ga _{\no}$ are empty.
\item D*D case $a_1=a_3=-1$: $R_{\no}$ is the set on and above the curve $\Ga _{\no}$ given by
\EQ{\label{GammanoDD}
\om = -\frac t3 + \frac 15 t^3, \quad \ga=-\frac 2{3t} - \frac 65 t,
\quad \sqrt{\frac53}<t <  \infty.
}
We have $(\om,\ga) \to (0,-\frac{8\sqrt{15}}{15})$ as $t \to \sqrt{\frac53}+$ and $(\om,\ga)\to (\infty,-\infty)$ as $t \to  \infty $.%
\end{enumerate}
Whenever $\Ga_\no$ exists (cases 1, 2, 4), it has a negative slope everywhere on the curve.
\end{theorem}

Note that, in the last  D*D case, $\frac{8\sqrt{15}}{15}\approx 2.06559$, which agrees with Figure \ref{log_region4}.

\begin{proof} By Lemma \ref{th3.2}, any $(\om,\ga)\in \Ga_\no$ must have a first positive zero $t$ which is a double zero.
By \S\ref{S:00}, $(\om,\ga)$ is given by \eqref{om-c.formula}. The last statement of negative slope follows from \eqref{negative-slope}.

\medskip

1. F*F case $a_1=a_3=1$. Since $G(0+)>0$ and $G(x)<0$ for $x\gg 1$, there is always a solution to $G(x)=0$, and (E1) is valid for every $(\om,\ga)$. When a positive zero $t$ of $G$ is a double zero,
 by \eqref{om-c.formula} and \eqref{om>0},
\[
\om = \frac t3 - \frac 15 t^3, \quad \ga=\frac 2{3t} + \frac 65 t,
\quad  x_0 = \frac{5}{6t}- \frac 12 t>0, 
\quad 0<t <  \sqrt{\frac53}.
\]
For (E2) to fail, we need $x_0\ge t$, i.e.~$t\le \frac{\sqrt{5}}3$. Thus $R_\no=\Ga_\no$ is parametrized by the above formula with $0 < t \le \frac{\sqrt{5}}3$.
In the interior of this interval, $dw/dt >0$ and $d\ga/dt <0$. 

\medskip

2. F*D case $a_1=1=-a_3$:  Since $G(0+)>0$ and $G(x)>0$ for $x\gg 1$,  for each fixed $\om>0$, the function $f(\ga)=\inf _{x>0}G(x)$ is continuous and nondecreasing in $\ga$ and at a minimal $\ga=\ga_0(\om)$ we have $f(\ga_0(\om))=0$. At this $\ga$, $G$ has a double zero. The curve $\ga=\ga_0(\om)$ describes $\Ga_0$.
The first positive zero of $G$ exists if and only if $\ga\le\ga_0(\om)$, i.e., $(\om,\ga)$ is on or below $\Ga_\no$.
When a positive zero $t$ of $G$ is a double zero,
 by \eqref{om-c.formula} and \eqref{om>0},
\[
\om = \frac t3 + \frac 15 t^3, \quad \ga=\frac 2{3t} - \frac 65 t,
\quad  x_0 = -\frac{5}{6t}- \frac 12 t<0, 
\quad 0<t <  \infty.
\]
As $x_0<0$, $t$ is always the first positive zero of $G$.
Thus $\Ga_\no$ is parametrized by the above formula with $0 < t <\infty$. In this interval, $dw/dt >0$ and $d\ga/dt <0$.

\medskip

3. D*F case $a_3=1=-a_1$: Since $G(0+)>0$ and $G(x)<0$ for $x\gg 1$, there is always a solution to $G(x)=0$, and (E1) is valid for every $(\om,\ga)$. When a positive zero $t$ of $G$ is a double zero,
 by \eqref{om-c.formula},
\[
\om = -\frac t3 - \frac 15 t^3, \quad \ga=-\frac 2{3t} + \frac 65 t.
\]
But there is no $t>0$ such that $\om(t)>0$. Thus (E2) is valid for every $(\om,\ga)$, and $R_\ex$ is the entire $\R_+\times \R$.

\medskip

4. D*D case $a_1=a_3=-1$. Since $G(0+)>0$ and $G(x)>0$ for $x\gg 1$,  for each fixed $\om>0$, the function $f(\ga)=\inf _{x>0}G(x)$ is continuous and nondecreasing in $\ga$ and at a minimal $\ga=\ga_0(\om)$ we have $f(\ga_0(\om))=0$. At this $\ga$, $G$ has a double zero. The curve $\ga=\ga_0(\om)$ describes $\Ga_0$.
The first positive zero of $G$ exists if and only if $\ga\le\ga_0(\om)$, i.e., $(\om,\ga)$ is on or below $\Ga_\no$.
When a positive zero $t$ of $G$ is a double zero,
 by \eqref{om-c.formula} and \eqref{om>0},
\[
\om = -\frac t3 + \frac 15 t^3, \quad \ga=-\frac 2{3t} - \frac 65 t,
\quad  x_0 = \frac{5}{6t}- \frac 12 t<0, 
\quad \sqrt{\frac53}<t <  \infty.
\]
As $x_0<0$, $t$ is always the first positive zero of $G$.
Thus $\Ga_\no$ is parametrized by the above formula with $\sqrt{\frac53} < t <\infty$. In this interval, $dw/dt >0$ and $d\ga/dt <0$. 
\end{proof}

Among the four cases, F*F is the most interesting, with standing waves on both sides of $\Ga_{\no}$.
Note that at the end point of $\Ga_\no$,
\EQ{\label{om1ga1}
(\om_1,\ga_1) =(\frac {2\sqrt5}{27},\frac{4\sqrt5}5) \approx (0.1656, 1.7889),\quad t =\frac{\sqrt{5}}3\approx 0.7454,
}
the first positive zero $t =\frac{\sqrt{5}}3$ is a triple zero of $G$.

\begin{lemma}
\label{th3.5}
In all cases  (F*F, F*D, D*D) when $\Ga_\no$ exists, 
for any $(\om_0,\ga_0)=(\om(t),\ga(t))$ given by \eqref{om-c.formula} on $\Ga_\no$, we have
\[
\lim_{\om\to \om_0-, \ga\to \ga_0-}  \phi_0(\om,\ga)=t.
\]
In the F*F case, we also have
\[
\lim_{\om\to \om_0+, \ga\to \ga_0+}  \phi_0(\om,\ga)=x_0
=  \frac{5}{6\, t} - \frac 12 t ,\quad x_0\ge t,
\]
with $x_0=t$ only at $(\om_1,\ga_1)$ with $t =\frac{\sqrt{5}}3$.
\end{lemma}
Note that the limits are taken through $(\om,\ga)\in R_\ex$ and
 exist by Lemma \ref{th3.1b}. %

\begin{figure}[H]
\centering
\includegraphics[width=0.47\textwidth, height=0.3\textwidth]{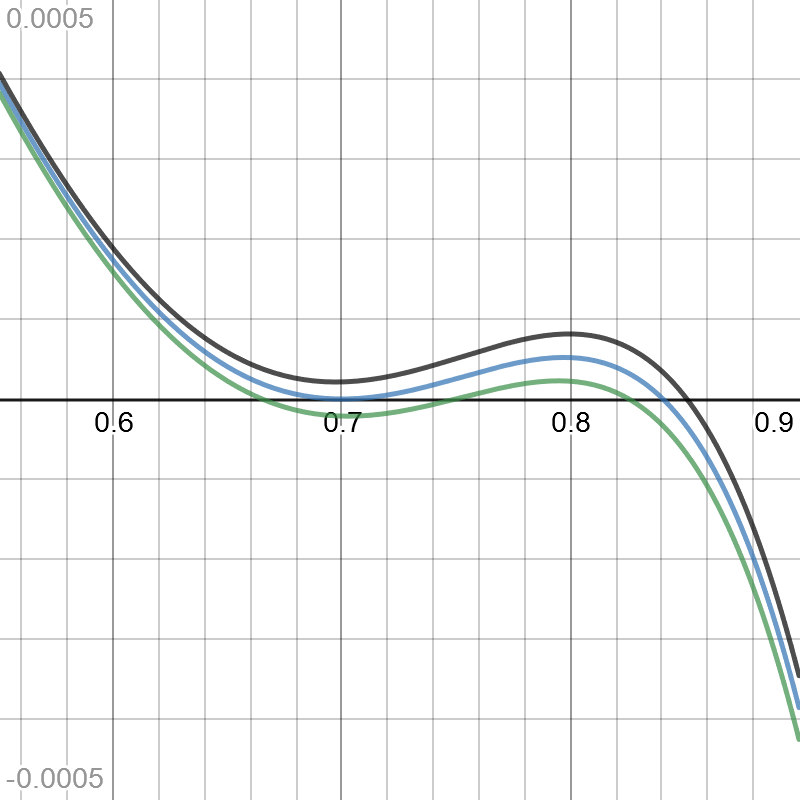}
\vskip2mm
\caption{Perturbations of $G_0(x)$ with double zero $t=0.7$, zero $x_0=0.8405$, FDF case.}
\label{G0-perturbations}
\end{figure}

\begin{proof}
By \eqref{G.factorization}, $G(x)=G_0(x)=-\frac{a_3}5x^2(x-t)^2(x-x_0)$ at $(\om_0,\ga_0)$. Thus for $(\om,\ga)$ near $(\om_0,\ga_0)$,
\[
G(x)=-\frac15x^2(x-t)^2(x-x_0)+\frac {\om-\om_0}2 x^2 + \frac {\ga-\ga_0}2 x^4.
\]

In both F*D and D*D cases, $a_3=-1$ and $x_0<0$. The  graph of $G_0$ is nonnegative for $0<x<\infty$ and its only positive zero is $x=t$. The zero near $t$ persists when $\om<\om_0$ and $\ga<\ga_0$, and converges to $t$ as $\om \to \om_0-$ and $\ga \to \ga_0-$. The zero disappears when $\om>\om_0$ and $\ga>\ga_0$.

In the F*F case, $a_3=1$ and $x_0\ge t$. For the non-endpoint case $t<x_0$,
the graph of $G_0$ is nonnegative for $0<x<x_0$ and it touches the $x$-axis at $x=t$ (see Figure \ref{G0-perturbations}).
 The zero near $t$ persists when $\om<\om_0$ and $\ga<\ga_0$, and converges to $t$ as $\om \to \om_0-$ and $\ga \to \ga_0-$.
The zero near $t$ disappears when $\om>\om_0$ and $\ga>\ga_0$, and the first positive zero jumps to near $x_0$ and converges to $x_0$ as $\om \to \om_0+$ and $\ga \to \ga_0+$. In the end point case $t=x_0$, the zero persists under any perturbation and converges to $t$ in the limits.
\end{proof}

\subsection{Solutions inside the standing wave homoclinic orbit}
\label{sec3.4}

In this Subsection we consider those $(\om,\ga)\in \R_+\times \R$ such that (E1) and (E2) holds, the first positive zero $\phi_0$ of $G$ exists and is not a double zero. We want to study solutions of the planar system \eqref{planar} that lie inside the homoclinic orbit passing $(\phi_0,0)$ on the phase plane.

The following lemma implies that,  in most cases, inside the homoclinic orbit passing $(\phi_0,0)$ on the phase plane, there is only one fixed point $(x_1,0)$, and any other solution is a periodic orbit around the fixed point. This is the familiar picture for a single power nonlinearity. However,  there are exceptions in the FDF case $a_1=a_3=1$ with $\ga>\sqrt 3$.

\begin{lemma}\label{th3.3}
Fix $\om>0$ and $\ga\in \R$.
If the first positive zero $\phi_0$ of $G$ exists and $g(\phi_0)<0$, there is a unique $x_1 \in (0,\phi_0)$ such that $g(x_1)=0$ in the F*D, D*F and D*D cases. It is also true in the F*F case if $-\I < \ga\le \sqrt 3$.
\end{lemma}

\begin{proof} A positive zero of $g(x)$ is also a  positive zero of 
\[
h(x) = \frac {g(x)}x = \om -a_1x+\ga x^2-a_3x^3.
\]
Note  $h(0)>0>h(\phi_0)$ and hence $h$ has a zero $x_1$ in $(0,\phi_0)$. If $h$ has more than one zero  in $(0,\phi_0)$, then
it has 3 zeros
\[
0<x_1 \le x_2 \le x_3 < \phi_0, \quad x_1 < x_3.
\]
Here we allow $x_1=x_2$ or $x_2=x_3$.
Then 
$h(x)$
must have a local minimum $x_-$ and a local maximum $x_+$  in
 $(0,\phi_0)$ with
\[
0< x_- < x_+ < \phi_0, \quad h(x_-)\le 0 \le h(x_+), \quad  h(x_-)< h(x_+).
\]
Thus $h'(x_-)=h'(x_+)=0$ and $x_\pm$ solve $3a_3x^2 -2\ga x +a_1=0$. Thus
\[
x_\pm = \frac1{3a_3}(\ga \pm \sqrt{\ga^2 -3a_1a_3}), \quad \ga^2 > 3a_1a_3.
\]
If $a_1a_3=-1$, then $x_+$ and $x_-$ have opposite sign and is invalid.  

We also have $h''(x_-) \ge 0 \ge h''(x_+)$. However, $h''(x)=2\ga-6a_3x$, and we must have $a_3=1$. This excludes the D*D case.

The only remaining case is F*F $a_1=a_3=1$. To avoid a contradiction to $x_\pm>0$, we need $\ga>\sqrt 3$.
\end{proof}

For the remaining FDF case $a_1=a_3=1$ and $\ga>\sqrt 3$, there may be exceptions.

\begin{figure}[H]
\centering
\includegraphics[width=0.47\textwidth, height=0.3\textwidth]{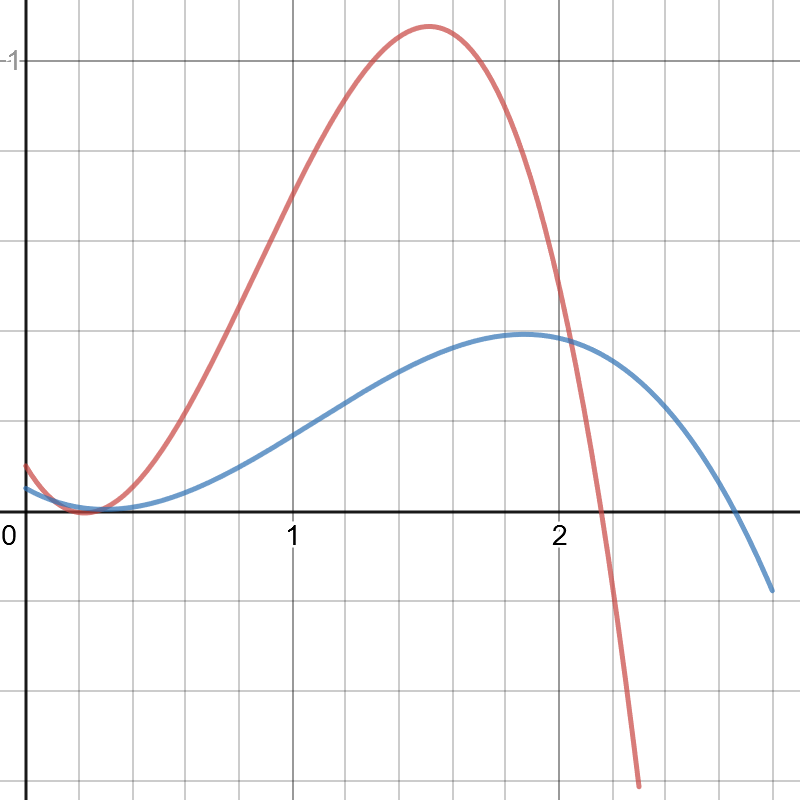}
\hfill
\includegraphics[width=0.47\textwidth, height=0.3\textwidth]{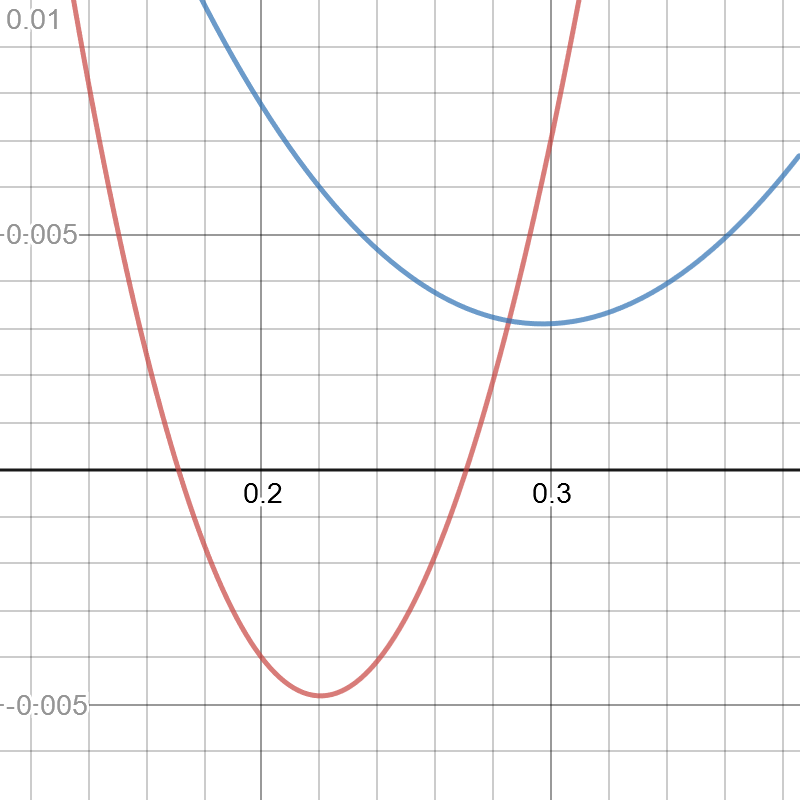}
\caption{$g(x)/x$ (red) and $G(x)/x^2$ (blue) for $0<x<2.8$ and $0.12<x<0.38$}
\label{ex3-5}
\end{figure}

\begin{example}\label{ex3.4}
Let $\om=0.1$ and $\ga=2.6$ in the FDF case. We have
\[
\frac{g(x)}x = 0.1 - x + 2.6x^2 -x^3, \quad 
\frac{G(x)}{x^2} = 0.05-\frac{x}{3}+\frac{2.6x^2}{4}-\frac{x^3}{5}.
\]
Numerically, the first positive zero of $G$ is $\phi_0 \approx 2.6585 $. The positive zeros of $g$ are approximately
\[
x_1=0.1711,\quad
x_2=0.2708,\quad
x_3=2.1581.
\]
See Figure \ref{ex3-5}.
In this example, there are two homoclinic orbits on the phase plane starting from $(x_2,0)$ and surrounding $(x_1,0)$ and $(x_3,0)$, respectively. Except these homoclinic orbits and fixed points, other solutions inside the standing wave orbit $(\phi, \dot \phi)$ are periodic orbits.\hfill \qed
\end{example}

\begin{example}\label{ex3.5}
A more systematic way is to consider
\[
g(x) = -x(x-x_1)(x-x_2)(x-x_3), \quad x_1=a, \quad x_2 = b, \quad x_3 = \frac {1-ab}{a+b},
\]
with parameters $a,b>0$, $ab<1$.
It is of the form \eqref{phi234.eq}, i.e., $g(x)=\om x - a_1 x^2 + \ga x^3 - a_3 x^4$, with $a_1=a_3=1$,
\[
\om=x_1x_2x_3= \frac {ab(1-ab)}{a+b}, \quad 
\ga = x_1+x_2+x_3=a+b+ \frac {1-ab}{a+b}.
\]
\begin{figure}[H]
\includegraphics[width=0.45\textwidth, height=0.3\textwidth]{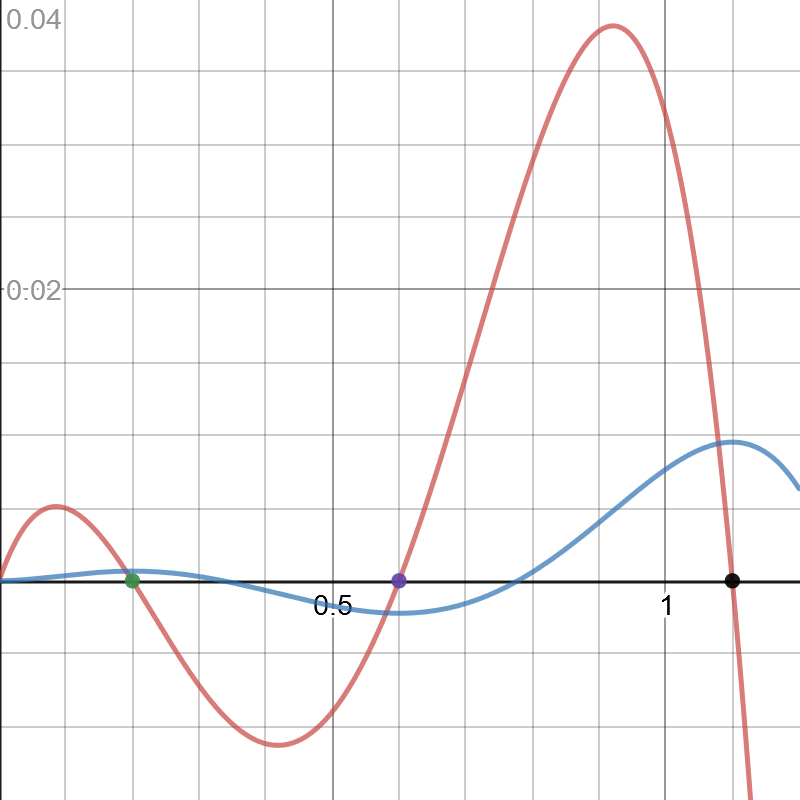}
(i) \hfill (ii)
\includegraphics[width=0.45\textwidth, height=0.3\textwidth]{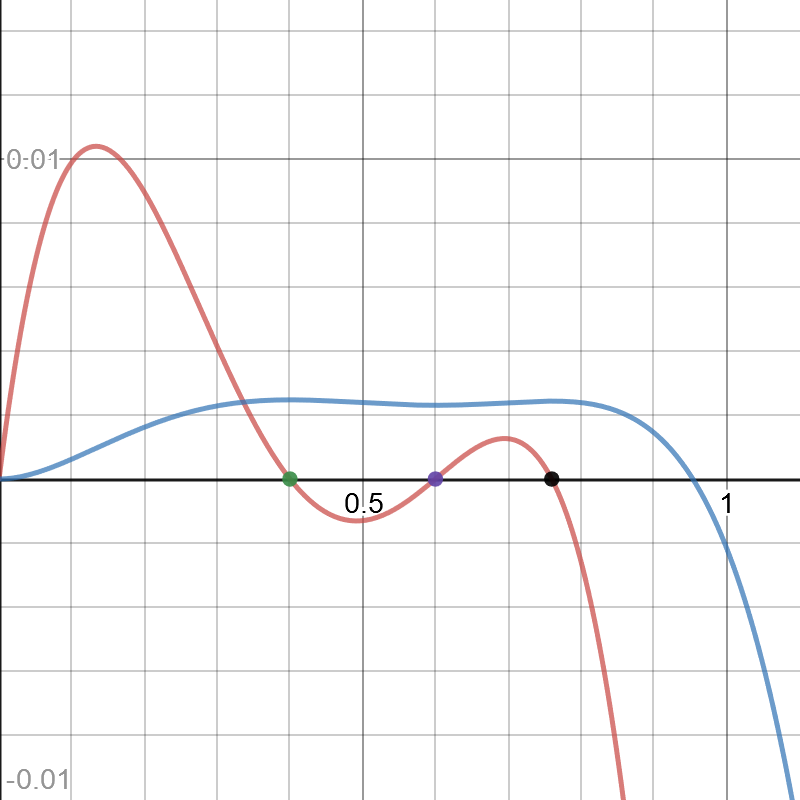}

\includegraphics[width=0.45\textwidth, height=0.3\textwidth]{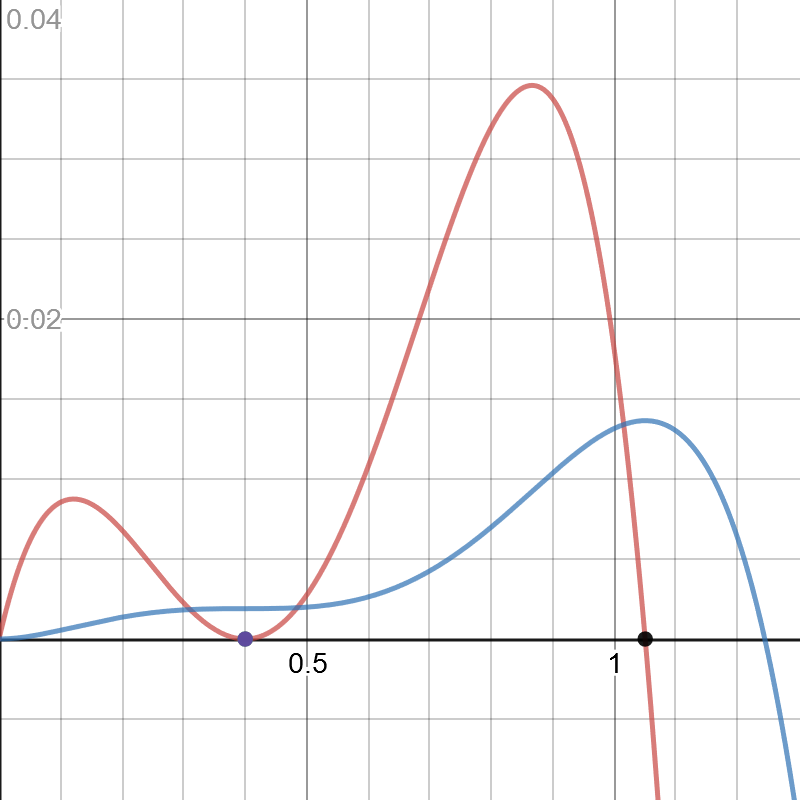}
{\small (iii)} \hfill {\small (iv)}
\includegraphics[width=0.45\textwidth, height=0.3\textwidth]{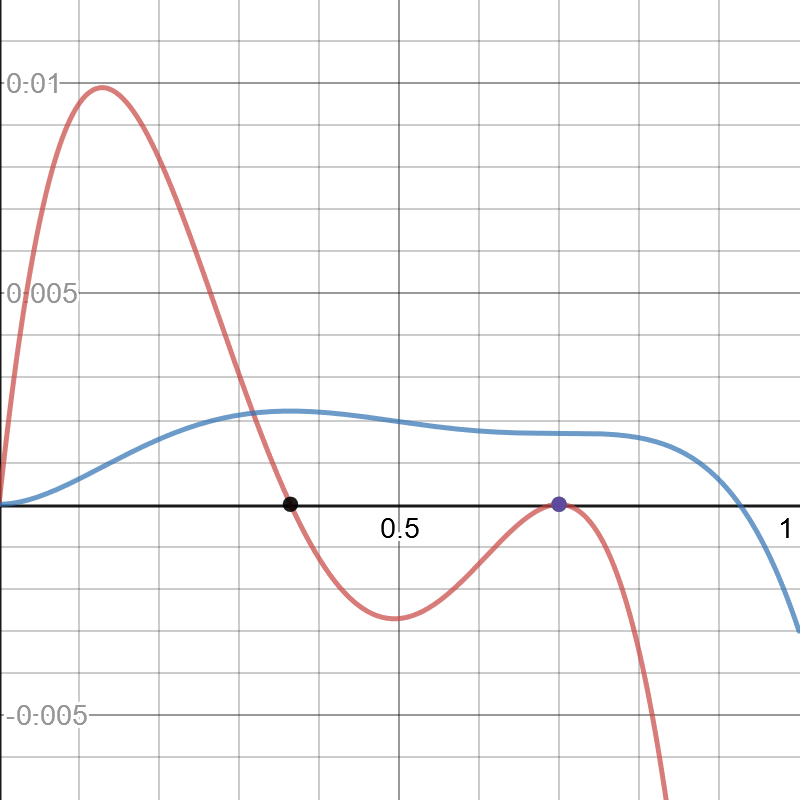}
\caption{$g(x)$ (red) and $G(x)$ (blue) for 
(i) $(a,b)=(0.2,0.6)$,
(ii) $(a,b)=(0.4,0.6)$, \newline
(iii)  $a=b=0.4$, 
(iv) $a=b=0.7$}
\label{ex3-6}
\end{figure}
For some choices of $a$ and $b$, $g$ has only one zero in $(0,\phi_0)$. See Figure \ref{ex3-6}(i) for the case $(a,b)=(0.2,0.6)$.
For some other choices of $a$ and $b$, $g$ has three zeros in $(0,\phi_0)$, and we may have $x_1=x_2$ or $x_2=x_3$. See Figure \ref{ex3-6}(ii,iii,iv). We can decouple the double zero of $g$ by decreasing $\om$ slightly in (iii) and by increasing $\om$ slightly in (iv).
If we fix $x_2$ and let $x_1 \to x_2-$, the homoclinic orbit around $(x_1,0)$ and all orbits inside it collapse into $(x_2,0)$. Similar phenomena happen if  we fix $x_2$ and let $x_3 \to x_2+$.\hfill \qed
\end{example}

In general, the Jacobian matrix of the planar system $\dot x=y$, $\dot y=g(x)$ is
\[
J = \mat{0 & 1 \\ g'(x) & 0}.
\]
At a fixed point $(x_1,0)$, the eigenvalues $\la$ satisfy $\la^2 = g'(x_1)$. Thus, inside $(0,\phi_0)$,
\EN{
\item if $g$ has exactly one positive zero $x_1$, then $g'(x_1)<0$, and $\la \in i \R$, and $(x_1,0)$ is a center. We expect peroidic orbits around it.
\item If $g$ has three distinct positive zeros $0<x_1<x_2<x_3<\phi_0$, then  $g'(x_1)<0$, $g'(x_2)>0$, and  $g'(x_3)<0$. Thus  $(x_1,0)$ and  $(x_3,0)$ are centers, while  $(x_2,0)$ is a saddle.
\item In the degenerate case $x_1=x_2$ or $x_2=x_3$, the phase protrait can be derived as a limit of 3-distinct zero case.
}

\section{Stability regions of parameters}
\label{sec4}

In this section we study the stability of standing waves for those $(\om,\ga)$ in the existence region $R_{\ex}$. In view of Proposition \ref{th1.2}, we define the \emph{stability functional}
\EQ{
J(\om,\ga) = \frac d{d\om}\int \phi_{\om,\ga}^2(x)\,dx, \quad (\om,\ga) \in R_{\ex}.
}
Equivalent integral formulas for $J$ are given in \eqref{th1.2-1}--\eqref{th1.2-1e}. We now divide the existence region $R_{\ex}$ to 3 subregions:
\begin{enumerate}
\item $\dis \Si_S = \bket{ (\om,\ga) \in R_{\ex}: J(\om,\ga)>0}$, stable region,
\item $\dis \Si_U = \bket{ (\om,\ga) \in R_{\ex}: J(\om,\ga)<0}$, unstable region,
\item $\dis \Ga_{cr} = \bket{ (\om,\ga) \in R_{\ex}: J(\om,\ga)=0}$, where stability may change.

\end{enumerate}

By Proposition \ref{th1.2}, $\Si_S$ corresponds to stable standing waves and $\Si_U$ corresponds to unstable standing waves. The set $\Ga_{cr}$ is 
the boundary between $\Si_S$ and $\Si_U$.  It is called the \emph{curve of stability change}, although in principle it may also contain points across which the standing waves do not change stability.

\subsection{Limits of the stability functional $J$ near the nonexistence curve}
\label{sec4.1}

\begin{proposition}\label{th4.1}
In all cases  (F*F, F*D, D*D) when $\Ga_\no$ exists, 
for any $(\om_0,\ga_0)=(\om(t),\ga(t))$ given by \eqref{om-c.formula} on $\Ga_\no$ except the end point $(\om_1,\ga_1)$ in the FDF case,
$J(\om_0,\ga_0)$ is undefined, and
\EQ{\label{th4.1-eq1}
\lim_{\om\to \om_0-, \ga\to \ga_0-}  J(\om,\ga)=+\infty.
}
In the F*F case, we also have
\EQ{\label{th4.1-eq2}
\lim_{\om\to \om_0+, \ga\to \ga_0+}  J(\om,\ga)=-\infty.
}
\end{proposition}
Note that the limits are taken through $(\om,\ga)\in R_\ex$.
Numerically, we also observe the above limits at $(\om_1,\ga_1)$ in the FDF case, and $\lim_{\om\to \om_1, \ga\to \ga_1}  J(\om,\ga)$ does not exist. See next Subsection.
We do not attempt to prove it as it is more technical due to the triple zero of $G$.

\begin{proof}
By \eqref{th1.2-1c}, 
\EQ{\label{th4.1-eq3}
J(\om,\ga) = \frac{-\sqrt 2}{4g(\phi_0)} \int_0^{\phi_0} \frac {I(x,\om,\ga)}{ \sqrt{G(x)}^3}  \,dx,
\quad 
}
where 
\EQN{
I(x,\om,\ga)&=6 \phi_0 x^2 G(x) + x^4 g(\phi_0) -\phi_0 x^3 g(x) .
}
Note that $\phi_0$, $g$ and $G$ all depend on $\om$ and $\ga$. The integral is improper and $I$ vanishes at both ends.

First consider the limit \eqref{th4.1-eq1}. By Lemma \ref{th3.5}, $\phi_0(\om_0-,\ga_0-)=\phi_0(\om_0,\ga_0)=t$. 
The integrand of \eqref{th4.1-eq3} is of order $O(1)$ for $x$ near 0, and of order $O((\phi_0-x)^{-1/2})$ for $x$ near $\phi_0$. However, the former is uniformly integrable but the latter is not since $G=O((t-x)^2)$ for $x$ near $\phi_0$ in the limit.
To study the integral for $x$ near $\phi_0$, denote
\EQN{
I_1(x,\om,\ga) & = \frac{I(x,\om,\ga)}{\phi_0-x}.
}
It is continuous in all variables. In the limits $\om\to \om_0-$ and $\ga\to \ga_0-$, it is close to 
\[
I_1(x,\om_0,\ga_0)= \frac{6 t x^2 G(x,\om_0,\ga_0) +0 -t x^3 g(x,\om_0,\ga_0)  }{t-x}.
\]
Using \eqref{G.factorization} that $G(x)=-\frac{a_3}5x^2(x-t)^2(x-x_0)$ at $(\om_0,\ga_0)$ and $g=G'$, we have
\EQN{
I_1(x,\om_0,\ga_0) &=  \frac{ t x^3  }{t-x} \frac{a_3}5\bke{ 2x^2(x-t)(x-x_0)+O((x-t)^2)}
\\
&=I_2(x)+O(x-t),
\\
I_2&=- \frac{2a_3}5  t  x^5(x-x_0),
}
for $x$ near $t$. In the F*F case, $a_3>0$ and $x_0\ge t$. In the **D case, $a_3<0$ and $x_0<0$. In all cases, there is $0<\de\ll1$ such that
\[
\inf_{t-\de<x<t} I_2(x)>0;\qquad
\frac12 I_2(x) < I_1(x,\om_0,\ga_0) < 2I_2(x), \quad (t-\de<x< t).
\]
By continuity, for all $\om<\om_0$ and $\ga<\ga_0$ sufficiently close to $(\om_0,\ga_0)$ so that $\phi_0<t$ is also close to $t$,
\[
\frac13 I_2(x) < I_1(x,\om,\ga) < 3I_2(x), \quad (t-\de<x< t).
\]
Decompose
\[
J(\om,\ga) = \frac{-\sqrt 2}{4g(\phi_0)} \bke{\int_0^{t-\de} + \int_{t-\de}^{\phi_0}}
\frac { (\phi_0-x)I_1}{ \sqrt{G}^3}  \,dx
 = \frac{-\sqrt 2}{4g(\phi_0)} \bke{ K_1(\om,\ga)+K_2(\om,\ga)}.
\]
$K_1(\om,\ga)$ is uniformly bounded in $(\om,\ga)$ while 
\[
K_2(\om,\ga)\ge \int_{t-\de}^{\phi_0(\om,\ga)} \frac { (\phi_0-x)I_2(x)}{3 \sqrt{G(x)}^3}  \,dx,
\]
which converges to $+\infty$ as $\om\to \om_0-$ and $\ga\to \ga_0-$, noting that $G(x)=O((x-t)^2)$ and the integrand is of order $O((x-t)^{-2})$ in the limits. Since $g(\phi_0)<0$, this shows \eqref{th4.1-eq1}.

\medskip

In the F*F case and $\om\to \om_0+$ and $\ga\to \ga_0+$, we have  $\phi_0(\om+,\ga+)=x_0$ by Lemma \ref{th3.5}. The integrand are bounded by $O(1)$ near $0$ and by $(\phi_0-x)^{-1/2}$ near $\phi_0$ and uniformly integrable. However, the integrand  is not uniformly integrable for $x$ near $t$ as $G$ has a double zero at $t$ in the limits: For $x$ near $t$, we have
\EQN{
I(x,\om,\ga) 
& \sim I(t,\om_0+,\ga_0+)=6 x_0 t^2 G(t) + t^4 g(x_0) -x_0 t^3 g(t)     = t^4 g(x_0)<0.
}
Above $G$ and $g$ are evaluated at $(\om_0,\ga_0)$.
Thus for $\de>0$ sufficiently small,
\[
 \int_{t-\de}^{t+\de} \frac {I(x,\om,\phi)}{ \sqrt{G(x)}^3}  \,dx \le \int_{t-\de}^{t+\de} \frac {t^4 g(\phi_0)}{2 \sqrt{G(x)}^3}  \,dx \to -\infty
\]
as $\om\to \om_0+$ and $\ga\to \ga_0+$. Thus the entire integral in \eqref{th4.1-eq3} converges to $-\infty$.  Since $g(\phi_0)<0$, this shows \eqref{th4.1-eq2}.
\end{proof}

\subsection{Level sets of the stability functional $J$ and stability regions}
\label{sec4.2}

In this Subsection we present our numerical results for the level sets of $J(\om,\ga)$ on the $\omega$-$\gamma$ half plane, for F*F, F*D, D*F and D*D cases. See
Figures \ref{FRF_levelsets}--\ref{DRD_levelsets}.

\subsubsection{Level sets for F*F case}\label{sec4.2.1}

By Theorem \ref{th3.4}, $R_\ex$ is the entire $\omega$-$\gamma$ half plane except the non-existence curve $\Ga_\no$ given by \eqref {GammanoFF}. 

\begin{figure}[H]
\begin{minipage}[b]{.5\textwidth}
\centering
\includegraphics[width=1\textwidth, height=1\textwidth]{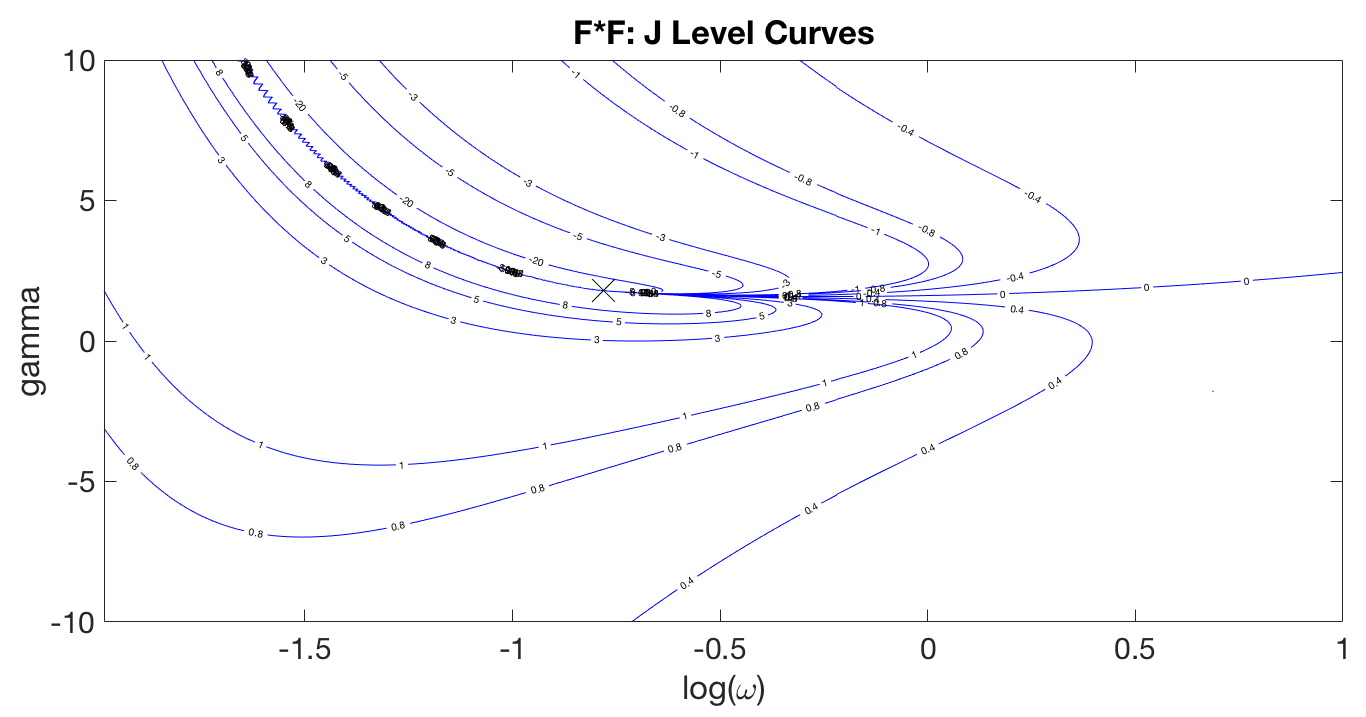}
\caption{(i) F*F case}
\label{FRF_levelsets}
\end{minipage}
\hfill
\begin{minipage}[b]{.5\textwidth}
\centering
\includegraphics[width=1\textwidth, height=1\textwidth]{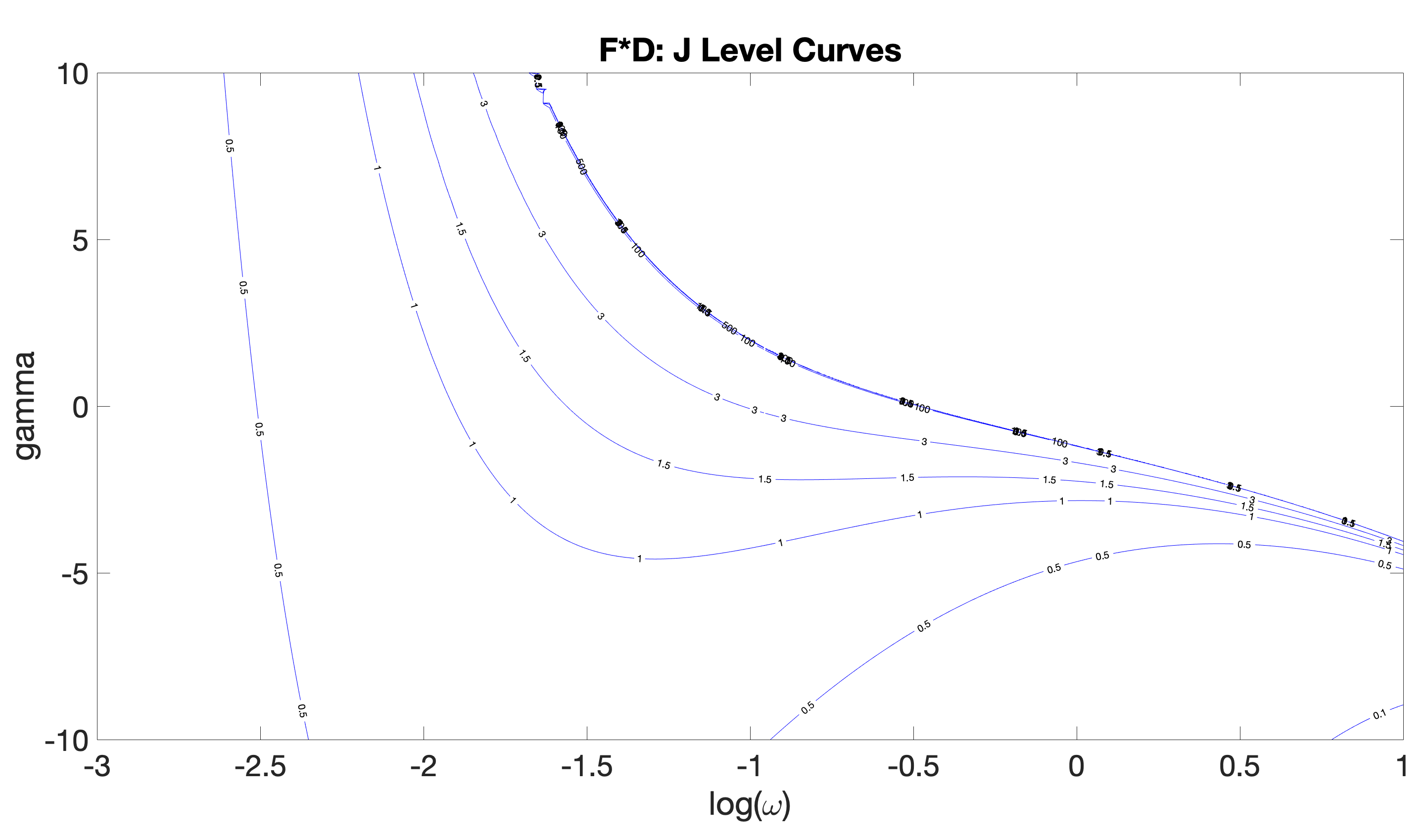} 
\caption{(ii) F*D case}
\label{FRD_levelsets}
\end{minipage}
\end{figure}

As shown in Figure \ref{FRF_levelsets}, also see Figure \ref{Fig2} in \S\ref{sec5.1},
the curve of stability change $\Ga_{cr}$, defined as the level set of $J=0$, is a graph of the form $\ga=\bar\ga(\om)$ defined for $\om_1<\om<\infty$. It emanates from the end point $(\om_1,\ga_1)\approx (0.1656, 1.7889)$ of $\Ga_\no$, given by \eqref{om1ga1}, and appears to have the same limiting slope given by \eqref{negative-slope},
\EQN{
\lim_{(\om,\ga)\in \Ga_\no \to (\om_1,\ga_1)}\frac{d\ga}{d\om} =\lim_{t\to \frac{\sqrt{5}}3} -\frac2{t^2}= - \frac{18}5.
}
The curve $\ga=\bar\ga(\om)$ appears to be smooth, and is decreasing until a critical point $(\om_2,\ga_2)$, whose numerical value is
\EQ{\label{rough-approx} 
(\om_2,\ga_2) \approx
(0.5548, 1.5817).
}
See \eqref{best-approx} of \S\ref{sec5.2} for a more accurate approximation.
The curve $\Ga_{cr}$ then becomes increasing for all $\om \in (\om_2, \infty)$. The value $\ga_2 \approx 1.5817$ is the global minimum of the function $\ga=\bar\ga(\om)$.

The union of $\Ga_\no\cup \Ga_{cr}$ appears to be a smooth and convex curve. Note that the smoothness and convexity of $\Ga_\no$, except at $(\om_1,\ga_1)$, follow from \eqref{negative-slope} and \eqref{GammanoFF}.
That of $\Ga_{cr}$ is only numerically observed.
For fixed  $\om \in (\om_1,\infty)$, $J(\om,\ga)$ is defined and continuous for all $\ga \in \R$.
It is possible to show that $J(\om,\ga_-)>0$ for some negative $\ga_-$, $J(\om,\ga_+)<0$ for some positive $\ga_+$, and hence $J(\om,\ga)=0$ for some $\ga \in (\ga_-,\ga_+)$ by the intermediate value theorem. However, it will require more work to show its uniqueness.

The stability functional $J(\om,\ga)$ has positive values below the union of $\Ga_\no\cup \Ga_{cr}$, and has negative values above it.  
Thus the region below the union of $\Ga_\no\cup \Ga_{cr}$ is the stable region $\Si_S$, and the region above it is the unstable region $\Si_U$.

For any fixed $c\not=0$, the level curve $J=c$ is a connected curve. It also comes out from $(\omega_1,\gamma_1)$, with the same limiting slope $ - \frac{18}5$. It diverges away from $\Gamma_{cr}$ as $\omega$ increases, with a negative slope. 

When $c>0$, the curve continues until it reaches to a specific $\omega$, and turns back toward $\omega = 0$, first with a positive slope, and then gradually switching to a negative slope. The value of $\gamma$ eventually goes to positive infinity, while $\om \to 0_+$.

When $c<0$, the curve continues and changes to a positive slope at some point, until it reaches to a specific $\omega$, and turns back toward $\omega = 0$ with a negative slope. The value of $\gamma$ eventually goes to positive infinity, while $\om \to 0_+$.

In both cases $c>0$ and $c<0$, the level curve $J=c$ emits from $(\omega_1,\gamma_1)$, changes slope sign twice, and eventually $(\om,\ga) \to (0_+,+\infty)$.

The picture agrees with Proposition \ref{th4.1} that, for any $(\om_0,\ga_0)\in \Ga_\no \setminus \{ (\om_1,\ga_1)\}$, we have
\[\lim_{\om\to \om_0-, \ga\to \ga_0-}  J(\om,\ga)=+\infty, \quad
\lim_{\om\to \om_0+, \ga\to \ga_0+}  J(\om,\ga)=-\infty.
\]

Since all level curves turn toward $\ga$-axis, it seems that the function $\om\mapsto \sup_{\ga \in \R} |J(\om,\ga)|$ is finite and decreasing in $\om$, with
\[
\lim_{\om \to \infty} \sup_{\ga \in \R} |J(\om,\ga)| = 0 .
\]

For stability, it follows from Proposition \ref{th1.2} that, for all $(\om,\ga)$ above $\ga = \bar \ga(\om)$, the solution $\phi(x)e^{i\om t}$ is orbital unstable. For all $(\om,\ga)$ below $\ga = \bar \ga(\om)$, the solution $\phi(x)e^{i\om t}$ is orbital stable. In particular, if
$-\I <  \ga < \ga_2$,  the solution $\phi(x)e^{i\om t}$ is orbital stable for all $\om>0$.
  
For the borderline standing waves $\phi_{\om,\bar \ga(\om)}e^{i\om t}$, $\om_1<\om<\infty$, by 
Comech and Pelinovsky \cite{MR1995870}, (also see 
\cite{CoCuPe} for gKdV), 
we expect them to be unstable for all $\om \not = \om_2$. It would be interesting to investigate the stability of the  standing wave $\phi_{\om_2,\ga_2}e^{i\om_2 t}$, since at $\ga=\ga_2$, all standing waves  $\phi_{\om,\ga_2}e^{i\om t}$ should be stable if $\om\not=\om_2$. The assumptions in \cite{MR1995870} for borderline standing waves likely fail for this degenerate case.

\subsubsection{Level sets for F*D case}

By Theorem \ref{th3.4}, $R_\ex$ is the region below the non-existence curve $\Ga_\no$ given by \eqref {GammanoFD}. 
Note that $R_{\ex}$ and $\Ga_{\no}$ are present in both FDD ($\ga>0$) and FFD ($\ga<0$) subcases.

As shown in Figure \ref{FRD_levelsets}, $J$ is positive in entire $R_\ex$. Hence every standing wave is orbital stable by Proposition \ref{th1.2}.

Agreeing with Proposition \ref{th4.1}, the values of $J$ converge to positive infinity as $(\om,\ga)$ converges to $\Ga_{\no}$.

\subsubsection{Level sets for D*F case}

By Theorem \ref{th3.4}, $R_\ex$ is the entire $\omega$-$\gamma$ half plane.

\begin{figure}[H]
\begin{minipage}[b]{.5\textwidth}
\centering
\includegraphics[width=1\textwidth, height=1\textwidth]{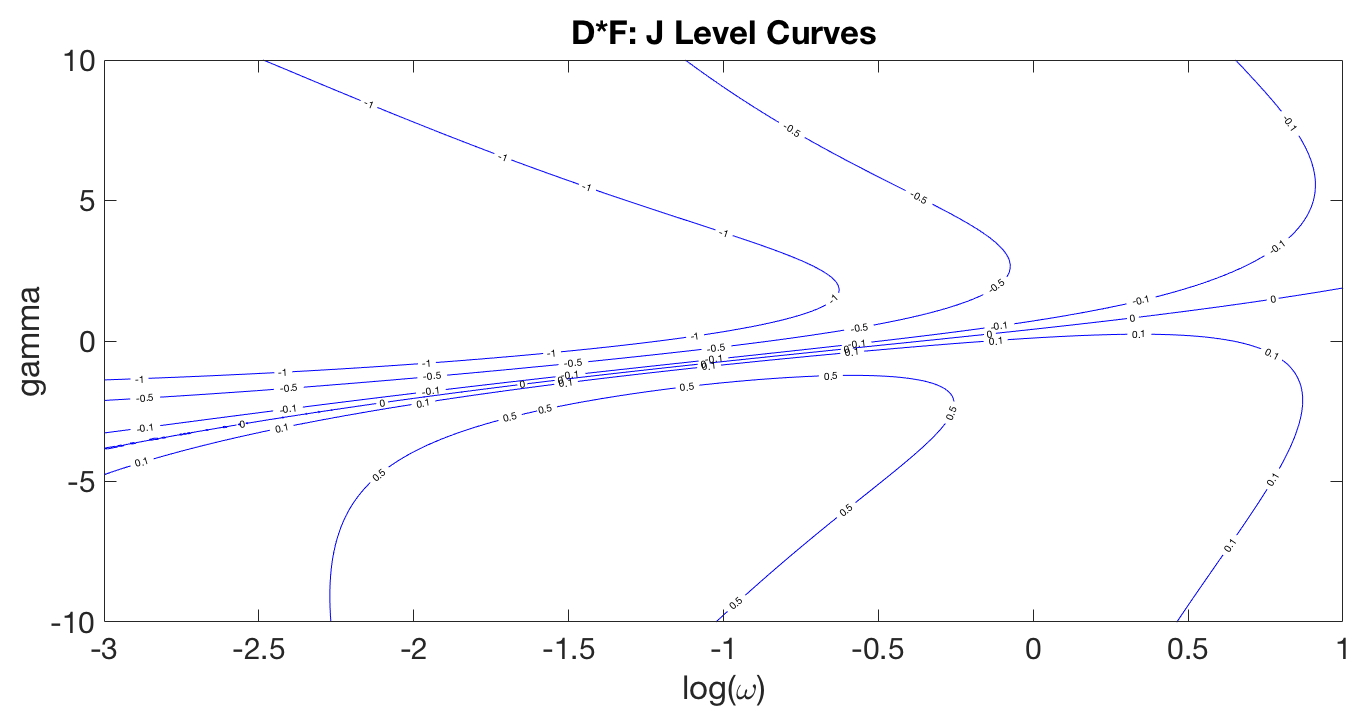}
\caption{(iii) D*F case}
\label{DRF_levelsets}
\end{minipage}
\hfill
\begin{minipage}[b]{.5\textwidth}
\centering
\includegraphics[width=1\textwidth, height=1\textwidth]{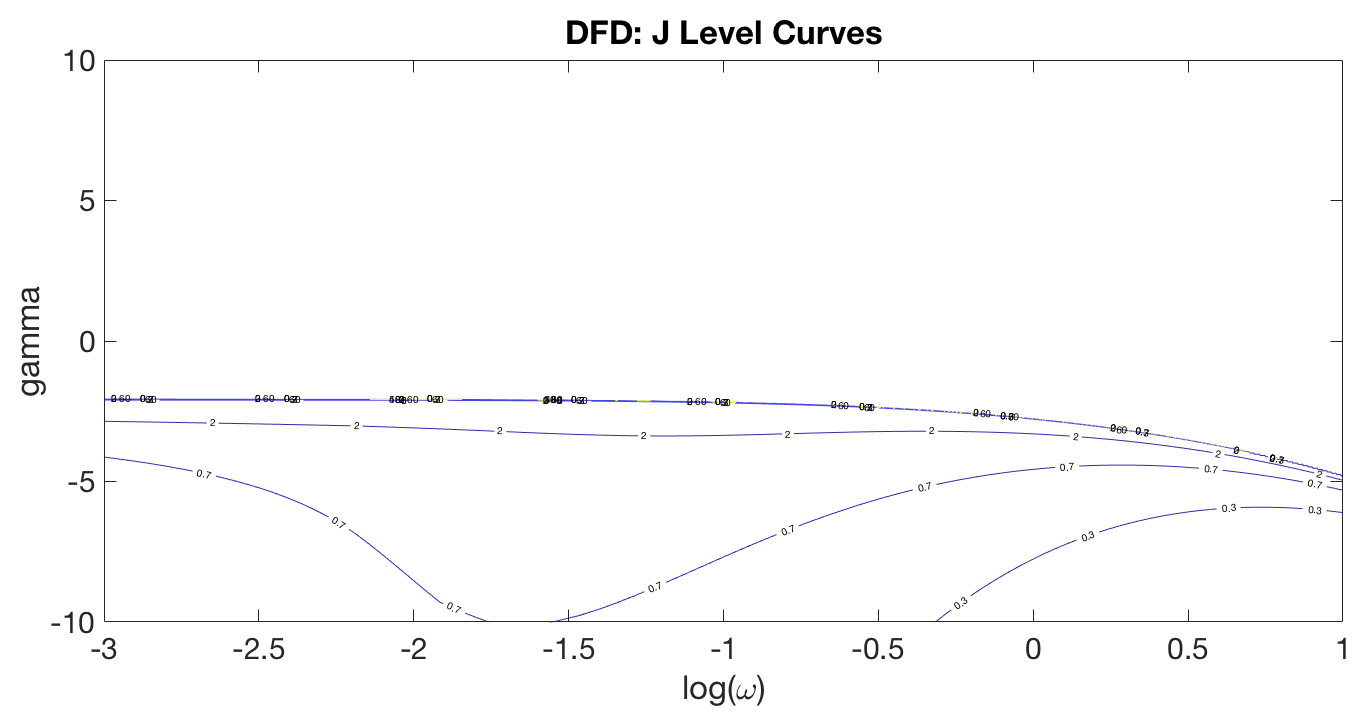} 
\caption{(iv) D*D case}
\label{DRD_levelsets}
\end{minipage}
\end{figure}

As shown in Figure \ref{DRF_levelsets}, the curve of stability change $\Ga_{cr}$, defined as the level set of $J=0$, is a graph of the form $\ga=\bar\ga(\om)$ defined for all $0<\om<\infty$. The curve is increasing for all $\om$. It appears to have a finite limit as $\om  \to 0_+$. Numerically, the closest point in our computation is
\[
(\om,\ga)=(0.001, -3.83459),
\]
as we cannot start from $\om=0$. The values of $J$ are positive below $\Ga_{cr}$, and negative above $\Ga_{cr}$. 
Thus the region below $\Ga_{cr}$ is the stable region $\Si_S$, and the region above $\Ga_{cr}$ is the unstable region $\Si_U$.

We expect the borderline standing waves $\phi_{\om,\bar \ga(\om)}e^{i\om t}$, $0<\om<\infty$ to be unstable by
Comech and Pelinovsky \cite{MR1995870},
if we could analyze $\ga=\bar\ga(\om)$ analytically.

The value of $J$ increases in magnitude as $(\om,\ga)$ diverges from the curve of stability change $\Ga_{cr}$.
Globally, we also observe
\[
\lim_{\om \to \infty} \sup_{\ga \in \R} |J(\om,\ga)| = 0 .
\]

\subsubsection{Level sets for D*D case}

By Theorem \ref{th3.4}, $R_\ex$ is the region below the non-existence curve $\Ga_\no$ given by \eqref {GammanoDD}. 
Note that $R_\ex$ and $\Ga_{\no}$ are present only in the DFD ($\ga<0$) subcase. Indeed, since $(\om,\ga) \to (0,-\frac{8\sqrt{15}}{15})$ as $t \to \sqrt{\frac53}+$, we have
\[
\ga < -\frac{8\sqrt{15}}{15}\approx -2.06559.
\]

As shown in Figure \ref{DRD_levelsets}, $J$ is positive in entire $R_\ex$. Hence every standing wave is orbital stable by Proposition \ref{th1.2}.

Agreeing with Proposition \ref{th4.1}, the values of $J$ converge to positive infinity as $(\om,\ga)$ converges to $\Ga_{\no}$.

\section{Numerics}
\label{sec5}

In this section we discuss our numerical methods and observations.

In Subsection \ref{sec5.1}, we describe the computations of the level curves of the stability functional $J(\omega,\gamma)$ and the stability regions for the four cases: F*F, F*D, D*F, and D*D in the right half $\omega$-$\gamma$ plane. The numerical results of all four cases agree with Theorem \ref{th3.4} and Proposition \ref{th4.1}. 

In Subsection \ref{sec5.2}, we describe the computation of the curve of stability change $\Gamma_{cr}$ and its minimal point $ (\omega_2, \gamma_2)$ in the F*F case. Our first method is based in successive zoom-in windows. Our second method is based on the \texttt{fsolve} command of MATLAB. Our most accurate approximation of $ (\omega_2, \gamma_2)$ is in \eqref{best-approx}.

In Subsection \ref{sec5.4}, we describe the 3 methods that we use for computing the standing wave $\phi_{\omega,\gamma}$. Our computed $\phi_{\omega,\gamma}$ agrees in the behaviour of the standing wave given by the planar dynamics. %

\subsection{Level curves and stability regions}
\label{sec5.1}
\begin{figure}[H]
\hspace{-25mm}
\includegraphics[width=1.25\textwidth, height=.65\textwidth]{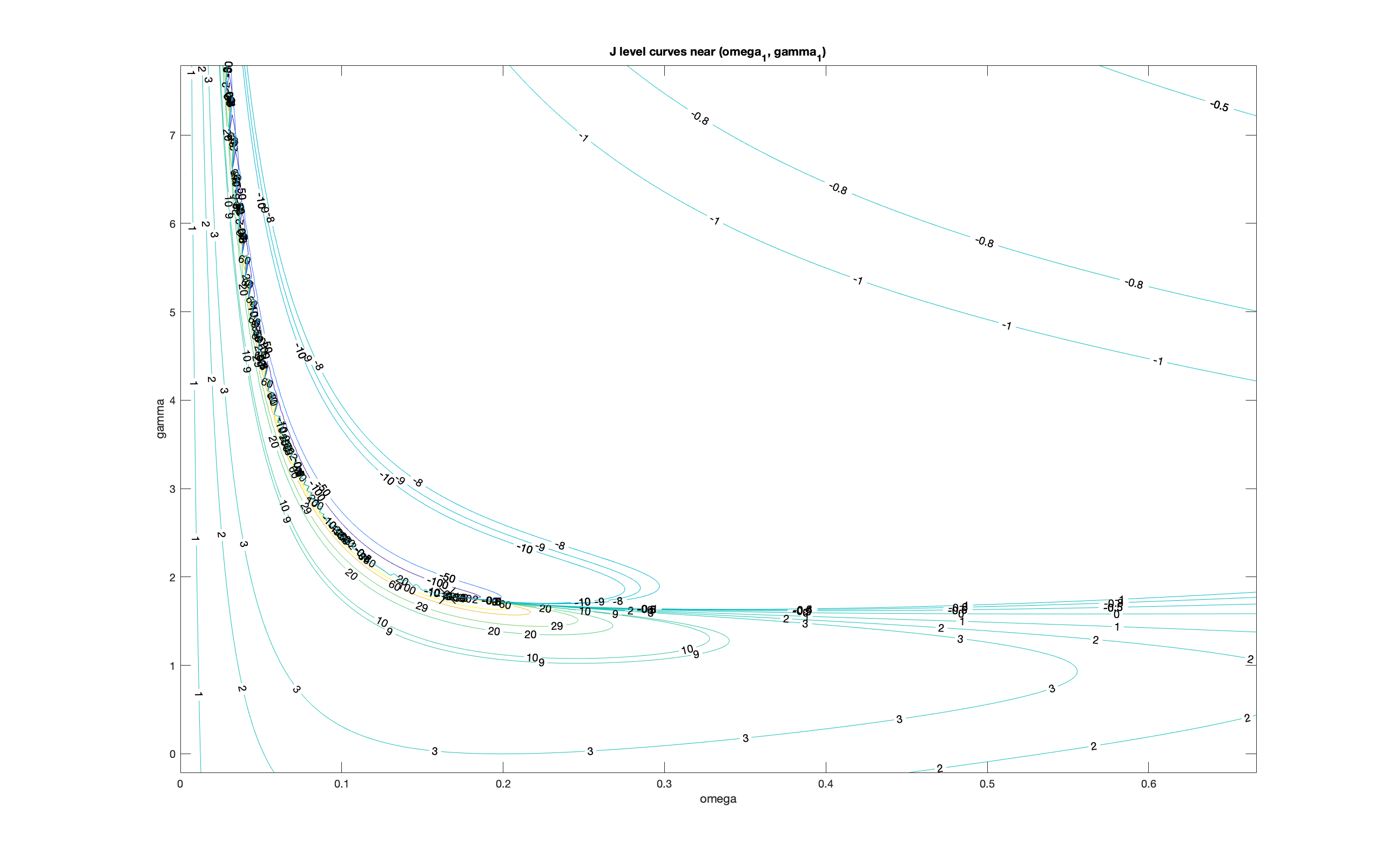}
\vspace{-10mm}
\caption{Level curves of $J$ near $(\omega_1,\gamma_1)$, F*F case.}
\label{Fig2}
\end{figure}

In this Subsection we describe how we obtain the level curves and stability regions numerically.
In MATLAB, we implement the existence condition (\ref{eq1.6}) into three ``if'' statements. The pairs $(\omega,\gamma)$ passing through all ``if'' statements are then categorized into stable, unstable, or boundary points by the value of $J(\omega,\gamma)$. We use (\ref{eq1.20.1}), the formula for $d''(\omega)$, to numerically compute $J$ in MATLAB.
The figures of the stability regions for all 4 cases are presented in Section \ref{sec4.2}.
Figure \ref{Fig2} is an additional plot of the level curves of $J$ in the F*F case.

\vspace{-4mm}

The level curves of $J$ are computed using MATLAB's \texttt{contour} function, and the mesh sizes for all cases are described below:

\medskip

For F*F case: For Figure \ref{FRF_levelsets}:

$(\omega,\gamma) \in [0.001,0.0104] \times [-10,10]$: $900 \times 900$ mesh points.

$(\omega,\gamma) \in [0.0104,0.603] \times [0,10]$: $900 \times 900$ mesh points, containing $(\om_1,\ga_1)$ and $(\om_2,\ga_2)$

$(\omega,\gamma) \in [0.0104,0.603] \times [-10,0]$: $900 \times 900$ mesh points.

$(\omega,\gamma) \in [0.603,5] \times [-10,10]$: $791 \times 900$ mesh points.

$(\omega,\gamma) \in [5,10] \times [0,10]$: $200 \times 200$ mesh points.

For Figure \ref{Fig2}:

$(\om,\ga) \in [0.001, 0.603] \times [0,10]$ : $500 \times 500$ mesh points.

\medskip

For F*D case: (Figures \ref{FRD_levelsets})

$(\omega,\gamma) \in [0.001,0.01] \times [-10,10]$: $450 \times 900$ mesh points.

$(\omega,\gamma) \in [0.01,0.0243] \times [0,10]$: $14 \times 900$ mesh points.

$(\omega,\gamma) \in [0.0243,0.01102] \times [0,10]$: $900 \times 900$ mesh points.

$(\omega,\gamma) \in [0.01102,1] \times [0,10]$: $807 \times 900$ mesh points.

$(\omega,\gamma) \in [0.01,1] \times [-10,0]$: $900 \times 900$ mesh points.

$(\omega,\gamma) \in [1,5] \times [0,10]$: $900\times 450$ mesh points.

$(\omega,\gamma) \in [1,10] \times [-10,0]$: $900\times 900$ mesh points.

$(\omega,\gamma) \in [5,10] \times [0,10]$: $200 \times 100$ mesh points.

\medskip

For D*F case: (Figures \ref{DRF_levelsets})

$(\omega,\gamma) \in [0.001,0.025] \times [-10,10]$: $900 \times 900$ mesh points.

$(\omega,\gamma) \in [0.025,5] \times [-10,10]$: $896 \times 900$ mesh points.

$(\omega,\gamma) \in [5,10] \times [-10,10]$: $200 \times 200$ mesh points.

\medskip

For DFD case: (Figures \ref{DRD_levelsets})

$(\omega,\gamma) \in [0.001,5] \times [-10,10]$: $900 \times 900$ mesh points.

$(\omega,\gamma) \in [5,10] \times [-10,10]$: $200 \times 200$ mesh points.

\medskip

When we used mesh points less than numbers stated above, the results were some non continuous level sets and zigzag curves,
 especially near  $\Gamma_{\no}$. However, using the refined mesh sizes, we get smooth looking connected level sets.

\subsection{The curve of stability change $\Gamma_{cr}$ and its minimal point $(\omega_2,\gamma_2)$}
\label{sec5.2}

Recall Subsection \ref{sec4.2.1} that, in the F*F case, the curve of stability change $\Gamma_{cr}$ is  numerically observed to be a graph $\ga = \bar \ga(\om)$, $\om_1<\om<\infty$, and its minimal $\ga$-value occurs at a critical point $(\omega_2, \gamma_2)$. We have
\[
(\om_1,\ga_1) =(\frac {2\sqrt5}{27},\frac{4\sqrt5}5) \approx (0.1656, 1.7889), \quad
(\om_2,\ga_2) \approx (0.5548, 1.5817).
\]
In this Subsection we describe how we get more accurate approximations of $\Gamma_{cr}$ and $(\omega_2, \gamma_2)$. 
Note that $\Ga_{cr}$ is only numerically observed without an analytic formula, and we have not proved its regularity, nor the existence of a unique minimum $\ga_2$ of $\bar \ga(\om)$.

For our computation, we use the MATLAB default precision of 16 decimal digits.

\medskip {\bf Method 1.}\quad
In our first method,
we improve our approximations by computation in a sequence of shrinking windows in the $\om$-$\ga$ parameter domain. See Table \ref{table1}.

\begin{table}[H]
\begin{center}
\begin{tabular}{ |c|c|c|c|c|c|c| } 
\hline
window		& $\om$ range & $\gamma$ range 		& mesh & $d\om$ & $d\ga$ & $\frac{d\om}{d\ga}$	 \\ 
\hline
$W_1$		& $[\omega_1+10^{-4}, 1.6656]$  	& $[1.55, 1.8]$		& $300 \times 1000$ & 0.005 & 2.5e-4 &20	 \\ 
\hline
$W_1'$		& $[0.4, 0.7]$  	&  $  [1.57, 1.6]$ 		& $60 \times 120 $ & 0.005 & 2.5e-4 &20	 \\ 
\hline
$W_2$		& $[0.5, 0.6]$  	& $[1.57,1.59]$		& $200 \times 800$ & 5e-4 & 2.5e-5 &20 \\ 
\hline
$W_2'$		& $[0.545, 0.565]$  	& $[1.5816,1.5818]$		& $40 \times 8$ & 5e-4 & 2.5e-5 &20 \\ 
\hline
$W_3$		& $[0.55,0.56]$  	& $ [1.58168,1.58172]$		& $400 \times 160$ & 2.5e-5 & 2.5e-7 &100	 \\ 
\hline
$W_3'$		& $[0.55,0.56]$  	& $[1.581704,1.581714]$		& $400 \times 40$ & 2.5e-5 & 2.5e-7 &100	 \\ 
\hline
\end{tabular}
\vspace{3mm}
\caption{Windows of computation. 2.5e-5 means $2.5\cdot 10^{-5}$.}\label{table1}
\end{center}
\end{table}

\subsubsection{Window $W_1$}
We start with the largest window 
\[
W_1 = [\omega_1 + 0.0001, 1.6656]\times [1.55, 1.8], \quad \text{mesh size: } 300 \times 1000.
\]
The choice of $W_1$ is based on the computation results for Figures \ref{FRF_levelsets} and \ref{Fig2}.
Since $\om_1 \approx 0.1656$,
\[
d\om\approx  0.005 ,\quad d\ga =0.000 25,\quad
\frac {d\om}{d\ga} = 20.
\]
The ratio $ {d\om}/{d\ga}=20$ is chosen larger than 1 in anticipation of a flat slope of $\ga=\bar \ga(\om)$ around $\om=\om_2$. 

For each mesh point $(\om_k,\ga_l)$ we compute $J(\om_k,\ga_l)$ using (\ref{eq1.20.1}).

\begin{figure}[H]
\hspace{-10mm}
\includegraphics[width=1.1\textwidth, height=.8\textwidth] {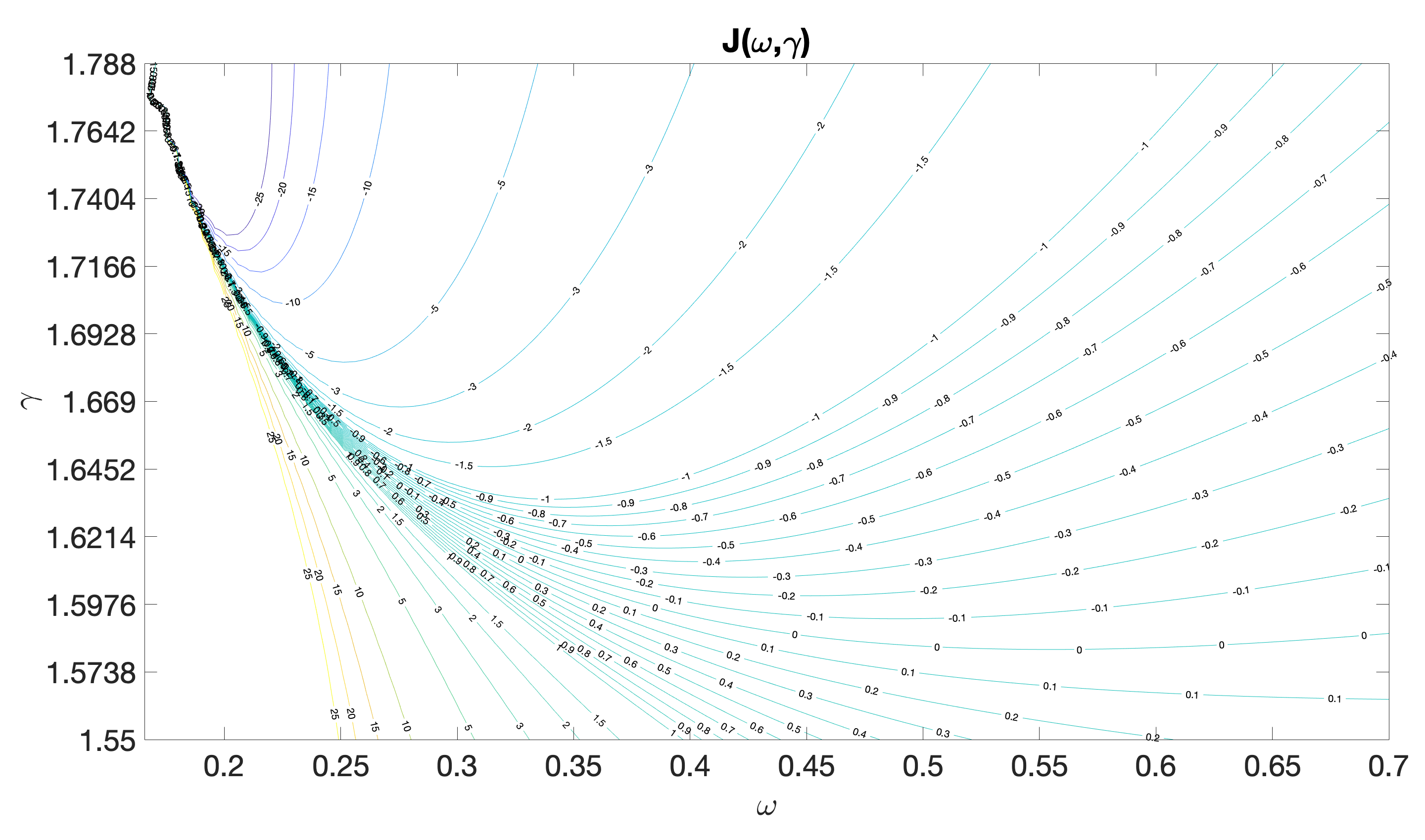}
\vspace{-6mm}  
\caption{Level curves of $J$ in the left half of window $W_1$, F*F case.}
\label{LevelsetsW1}
\end{figure}

In view of Figures \ref{FRF_levelsets} and \ref{Fig2}, we anticipate that $J(\om_k,\ga_l)$ is decreasing in $l$ at most points, except those close to $(\om_1,\ga_1)$ where some level curves of $J$ turn around. This is reinforced by Figure \ref{LevelsetsW1} on the level curves of $J$ in $W_1$. 
One observes that the portion of the level curve $J=c$ after turn around does not intersect $W_1$ if $|c|<20$. Hence $J(\om_k,\ga_l)$ is decreasing in $l$ in $W_1$ as long as $|J(\om_k,\ga_l)|<20$ (which is always the case for $\om\ge0.4$).

\begin{figure}[H]
\hspace{-3mm}
\includegraphics[width=1\textwidth, height=.45\textheight]{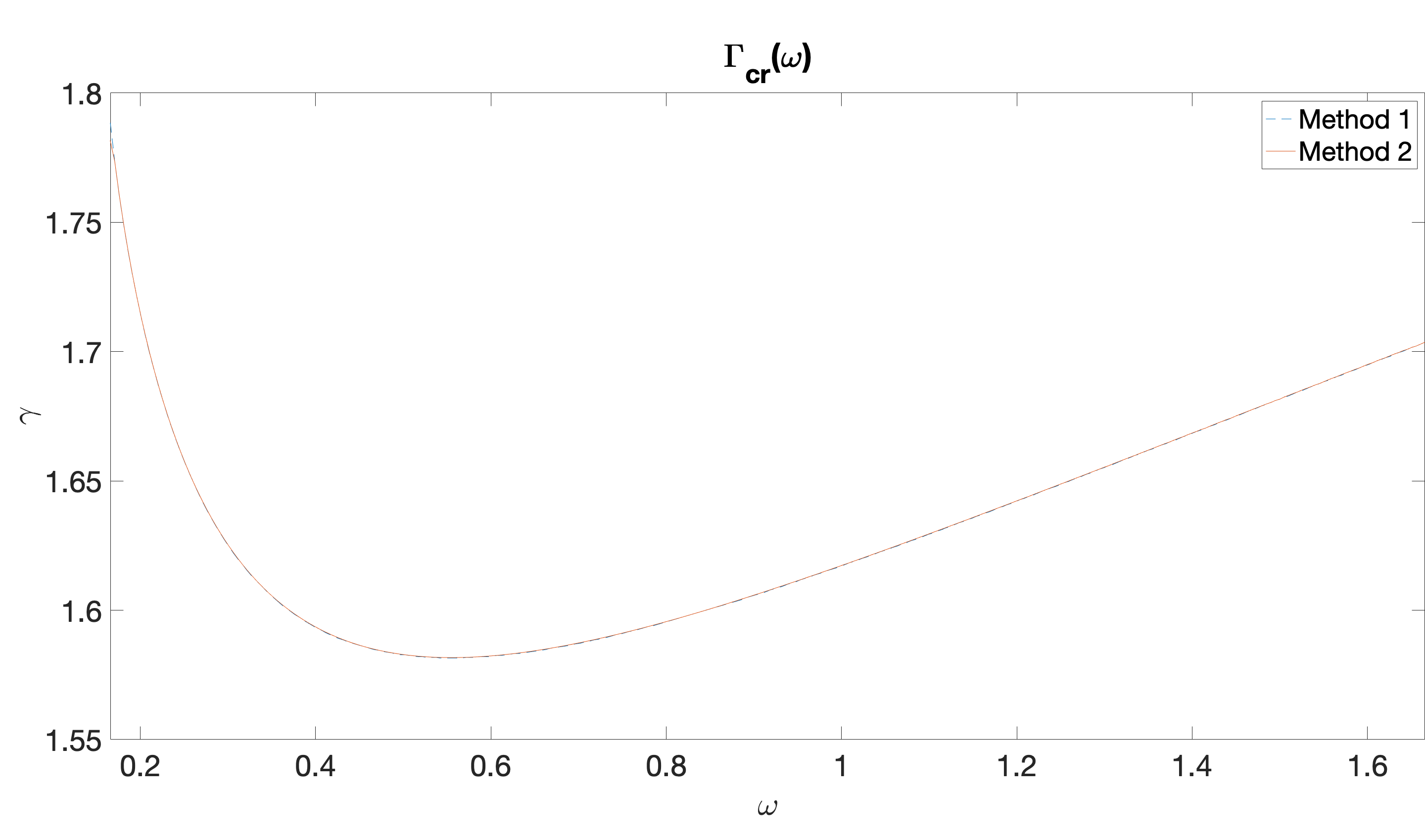}
\vspace{-1mm}
\caption{$\Gamma_{cr}$ in $W_1$.}
\label{Fig_improved_connection_no_cr}
\end{figure}

Thus, for each $\omega_k$ we define $\bar {\gamma}(\omega_k)= \ga_l$ where $\ga_l$ is the smallest such that
\EQ{
J(\omega_k,\gamma_l) \geq 0 > J(\omega_k,\gamma_{l+1}).
}
This gives us a discrete curve $(\om_k, \bar \ga(\om_k))_k$, shown as the dotted curve in Figure \ref{Fig_improved_connection_no_cr}.

We then define the approximate $\ga_2$ as
\EQ{
\gamma_2^*=\min_k \{\bar \gamma(\om_k)\},
}
and the $\om$-interval that attains $\ga_2^*$
\EQ{
[\om^*_{\min} , \om^*_{\max}] =  \{ \om_k: \ \bar \gamma(\om_k) = \gamma_2^*\}.
}
We expect our true values of $(\om_2,\ga_2)$ satisfy
\[
\bar \ga \le \ga_2 \le \bar \ga+d\ga, \quad \om^*_{\min}-d\om \le \om_2 \le \om^*_{\max}+d\om.
\]
We may take $\om_2^* =\frac12(\om^*_{\min}+ \om^*_{\max})$ as an approximation of $\om_2$. We obtain%
\EQ{
W_1: \quad \gamma_2^* =1.581531531531532,\quad
\omega_2^* \in [0.541955066755841,0.567036426872898].
}

\begin{figure}[H]
\hspace{-3mm}
\includegraphics[width=1.05\textwidth, height=.45\textheight]{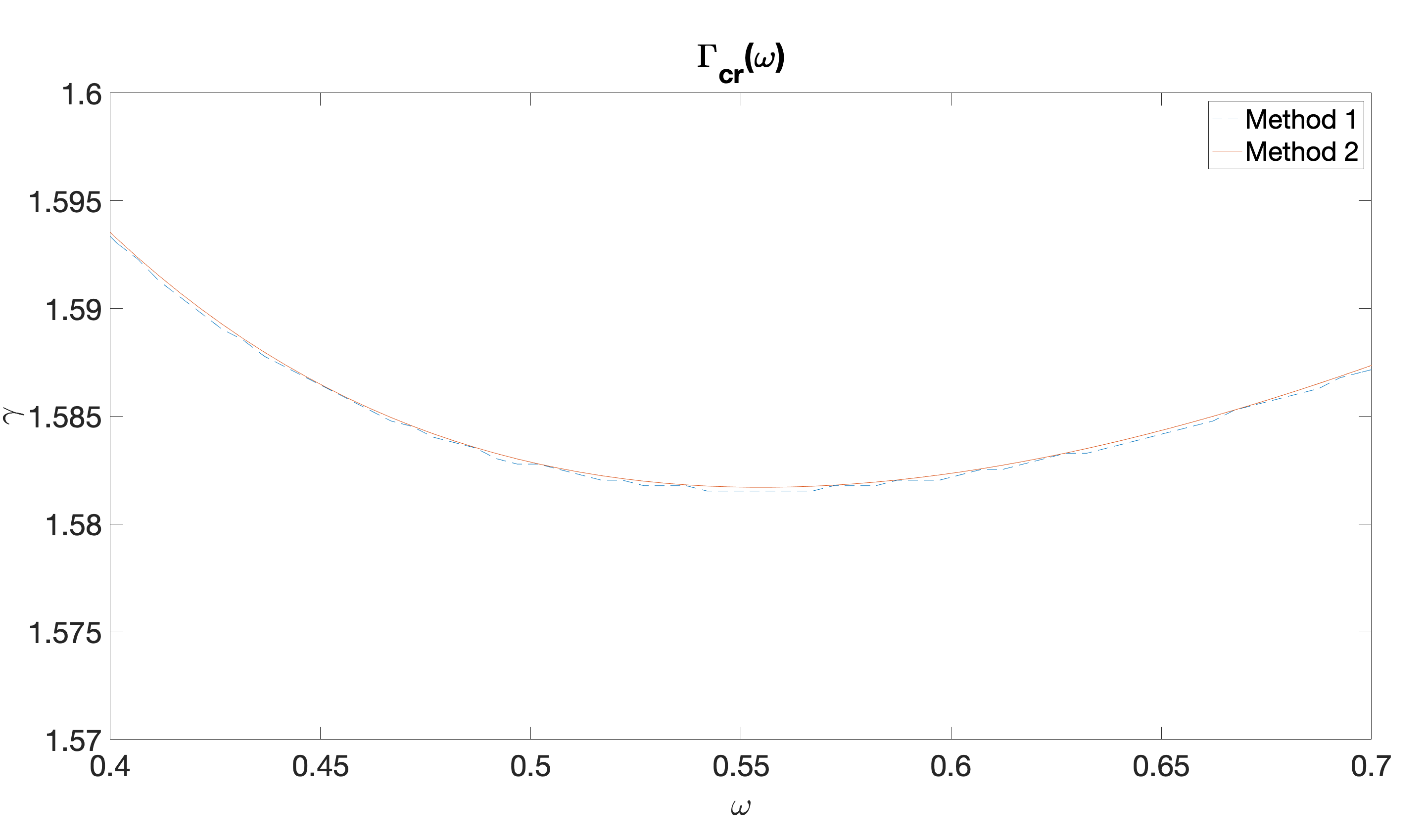}
\vspace{-4mm}
\caption{$\Gamma_{cr}$ in $W_1'$ }
\label{W1p_improved_connection_no_cr}
\end{figure}

From Figure \ref{Fig_improved_connection_no_cr}, we choose a subwindow $W_1'$ that clearly contains the line segment $[\om^*_{\min} , \om^*_{\max}] \times \{ \ga_2^*\}$ and graph a zoomed picture in $W_1'$ (using the previous computation result in $W_1$), see Figure \ref{W1p_improved_connection_no_cr}. Using this finer picture we will choose our next window $W_2$.

\bigskip
 
{\bf Method 2}. \quad Our second method is based on  MATLAB's \texttt{fsolve} command, which solves a zero of a given function near a given initial point. It can be considered as a refinement of Method 1 since we use the results of Method 1 as our initial points. However, there are many other ways to choose initial points. For example, with an approximation curve $(\bar \om_k,\bar \ga_k)$ from a previous computation, we can insert $N$ points between $(\bar \om_k,\bar \ga_k)$ and $(\bar \om_{k+1},\bar \ga_{k+1})$ by linear interpolation for each $k$, and use the portion of this new curve that is near the minimal point as our new initial points.

Specifically,
we use the 300 points
\[
(\om_k,\ga(\om_k)), \quad 1 \le k \le 300
\]
obtained in Method 1
as the initial points, and use  MATLAB's \texttt{fsolve} command with \texttt{levenberg-marquardt} algorithm option, to solve 
300 zeros
\[
(\bar \om_k,\bar \ga_k), \quad 1 \le k \le 300
\]
of the stability functional $J(\omega,\gamma)$ that are close to $(\om_k,\ga(\om_k))$. This parametric, discrete curve is shown as the solid curve in Figure \ref{Fig_improved_connection_no_cr}. It is a finer approximation of the curve $\ga=\bar \ga(\om)$ since the default precision of MATLAB when applying \texttt{fsolve} is 16 decimal digits, much less than our $d\ga$. These points $(\bar \om_k,\bar \ga_k)$ are not mesh points, but it is fine.

We now define the approximate $\ga_2$ and $\om_2$ as
\EQ{
\gamma_2^{**}=\min\nolimits_k \{\bar \gamma_k\},\quad
\om_2^{**}=\text{argmin}_k \{\bar \gamma_k\}.
}
This time $\om_2^{**}$ is attained at a single value of $\om$, not an interval as in Method 1.
We obtain
\EQ{\label{W1-approx}
W_1: \quad (\omega_2^{**}, \gamma_2^{**}) = (0.557003765501051,1.58170639609768).
}

\vspace{-5mm}  

\subsubsection{Window $W_2$}

Based on Figure \ref{W1p_improved_connection_no_cr} on the curve $\Ga_{cr}$ in window $W_1'$, we choose our second window
\[
W_2 = [0.5, 0.6]\times [1.57,1.59], \quad \text{mesh size: } 200 \times 800.
\]

We do not increase the ratio $d\om/d\ga=20$ since the mesh of $W_1'$ in Table \ref{table1} has more vertical points, which indicates that the ratio 20 is sufficient.

\begin{figure}[H]
\hspace{-3mm}
\includegraphics[width=1\textwidth, height=.45\textheight]{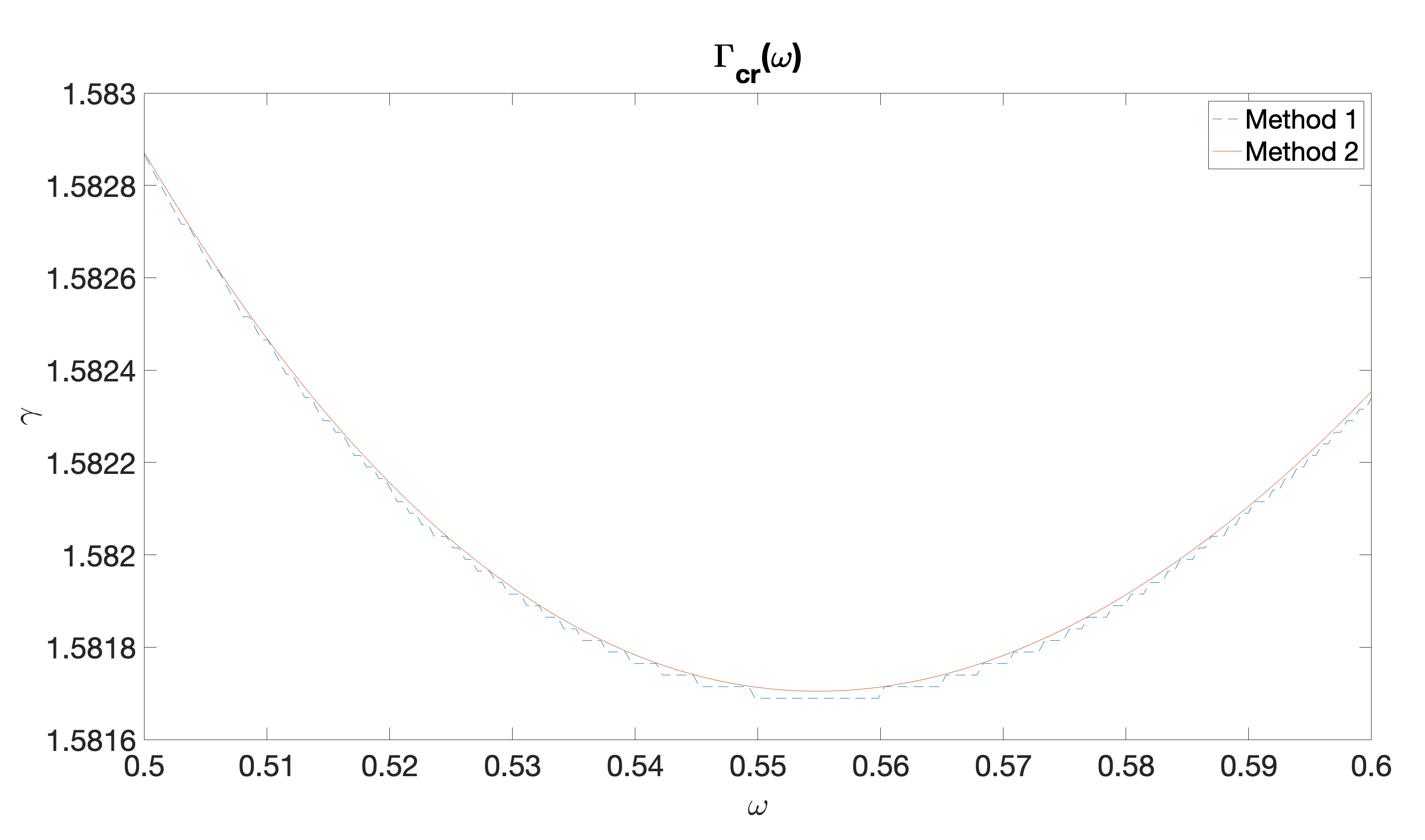}
\vspace{0mm}
\caption{$\Gamma_{cr}$ in $W_2$ }
\label{W2-Gcr}
\end{figure}

\begin{figure}[H]
\hspace{-3mm}
\includegraphics[width=1\textwidth, height=.45\textheight]{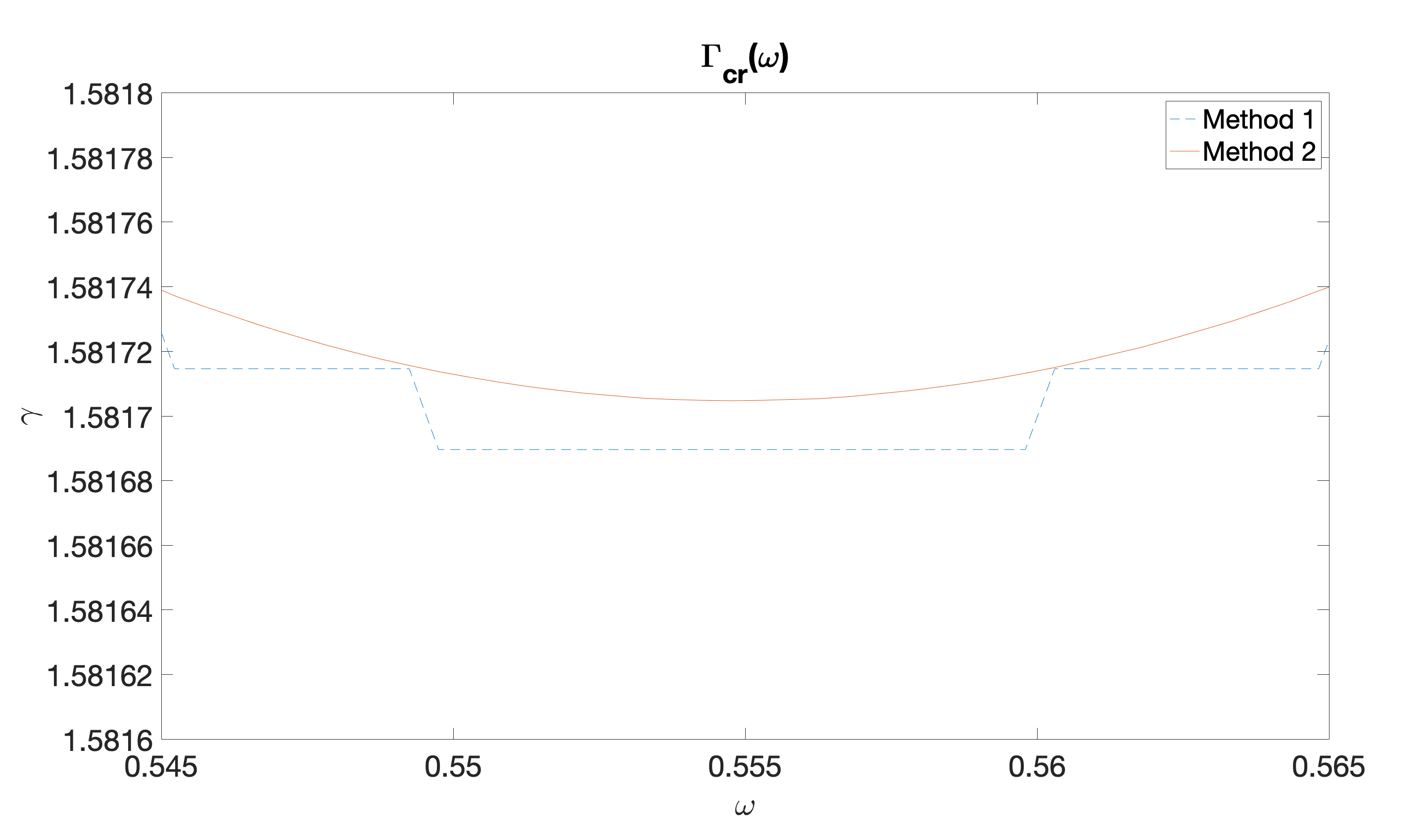}
\vspace{0mm}
\caption{$\Gamma_{cr}$ in $W_2'$ }
\label{W2p_improved}
\end{figure}

We repeat Methods 1 and 2 for window $W_2$. Note that $W_2$ is away from $(\om_1,\ga_1)$ and $J(\om_k,\ga_l)$ is decreasing in $l$ at \emph{all} mesh points in $W_2$. The approximation curves based on Methods 1 and 2 are presented in Figure \ref{W2-Gcr}. We also choose a subwindow $W_2'= [0.545, 0.565]\times [1.5816,1.5818]$, presented in Figure \ref{W2p_improved}, for the choice of $W_3$.

We obtain, by Method 1:
\EQ{
W_2: \quad \gamma_2^* =1.58168961201502,\quad
\omega_2^* \in [0.549748743718593 ,0.559798994974874],
}
and by Method 2:
\EQ{\label{W2-approx}
W_2: \quad (\omega_2^{**}, \gamma_2^{**}) = (0.554773875083001,1.58170476081013).
}

\subsubsection{Window $W_3$}

Based on Figure \ref{W2p_improved} on the curve $\Ga_{cr}$ in window $W_2'$, we choose our third window
\[
W_3 = [0.55,0.56]\times  [1.58168,1.58172], \quad \text{mesh size: } 400 \times 160.
\]

We increase the ratio $d\om/d\ga$ from $20$ to $100$ since the mesh of $W_2'$ in Table \ref{table1} has 5 times horizontal points than vertical points.

\begin{figure}[H]
\hspace{-3mm}
\includegraphics[width=1\textwidth, height=.45\textheight]{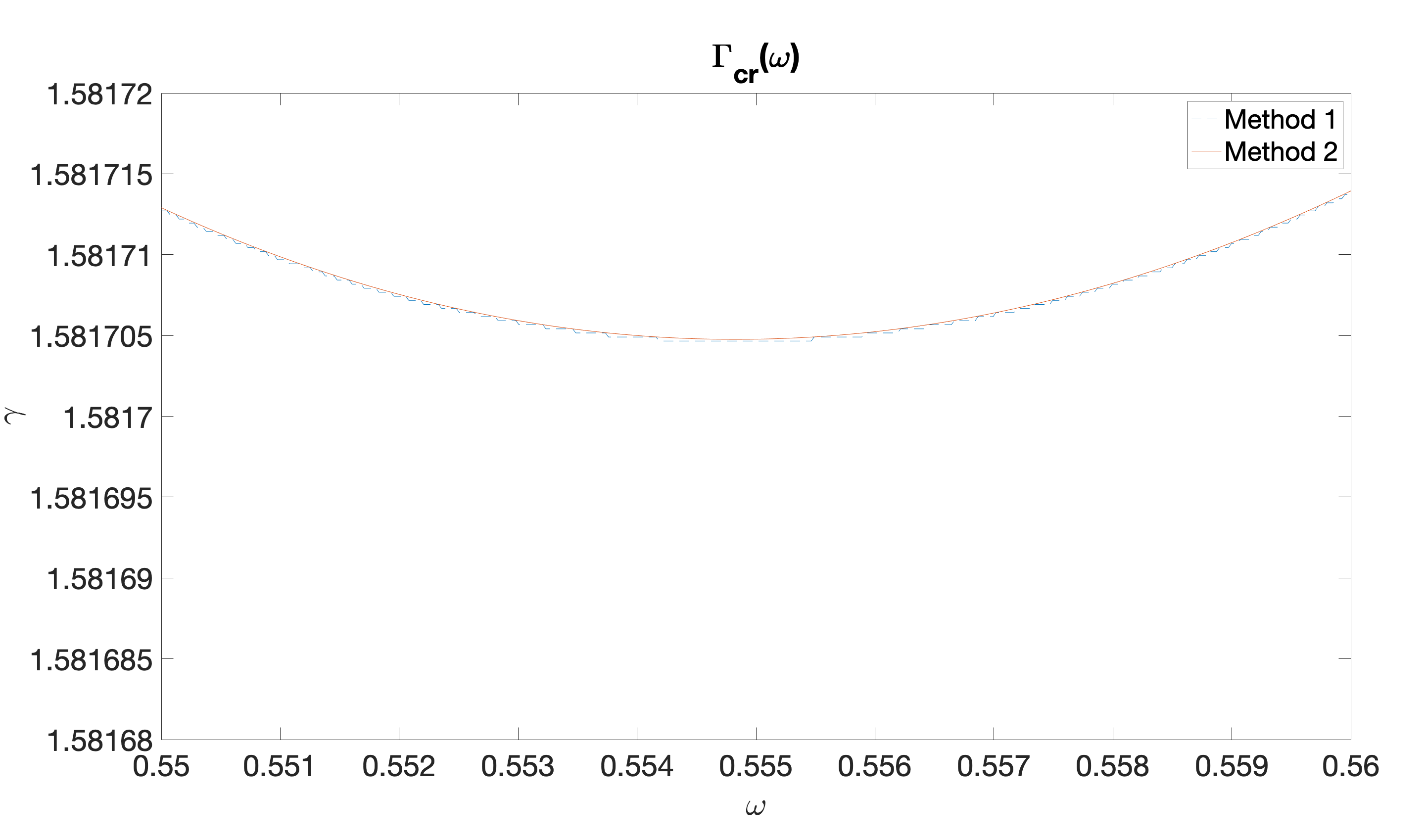}
\vspace{0mm}
\caption{$\Gamma_{cr}$ in $W_3$ }
\label{W3-Gcr}
\end{figure}

\begin{figure}[H]
\hspace{-3mm}
\includegraphics[width=1\textwidth, height=.45\textheight]{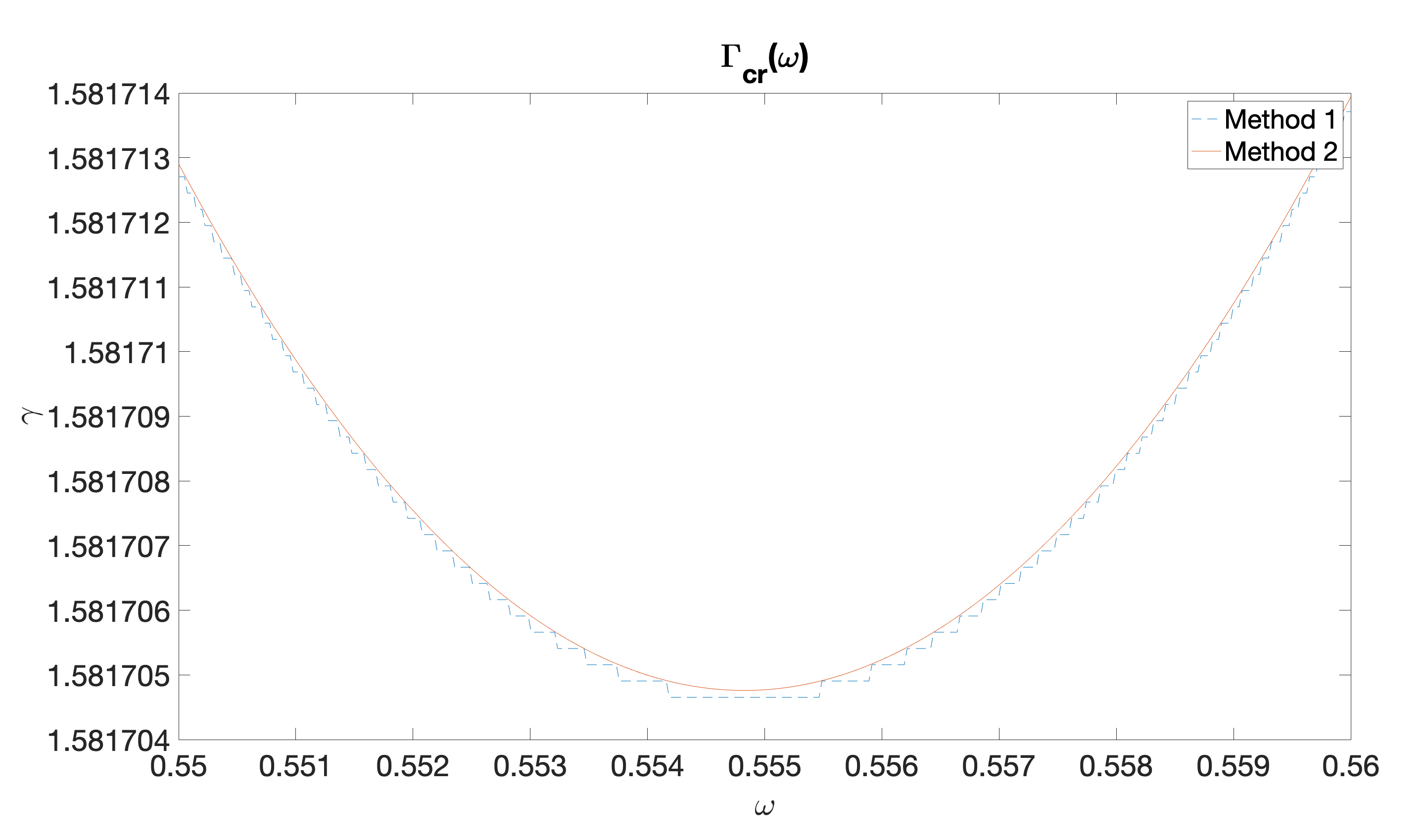}
\vspace{0mm}
\caption{$\Gamma_{cr}$ in $W_3'$ }
\label{W3p_improved}
\end{figure}

We repeat Methods 1 and 2 for window $W_3$. The approximation curves based on Methods 1 and 2 are presented in Figure \ref{W3-Gcr}. We also choose a subwindow $W_3'= [0.545, 0.565]\times [1.581704,1.581714]$, presented in Figure \ref{W3p_improved}, for a better local view.

We obtain, by Method 1:
\EQ{
W_3: \quad \gamma_2^* =1.58170465408805,\quad
\omega_2^* \in[ 0.554185463659148,0.55546365914787],
}
and by Method 2:
\EQ{\label{best-approx}
W_3: \quad (\omega_2^{**}, \gamma_2^{**}) =  (0.554837092755109,1.58170475989899).
}
Eq.~\eqref{best-approx} is our best approximation of $(\om_2,\ga_2)$. We collect $(\omega_2^{**}, \gamma_2^{**}) $ in 3 windows in Table \ref{table2}. Also included are $\De \omega_2^{**}$ and $\De \gamma_2^{**}$, the difference of the current values of $\omega_2^{**},\gamma_2^{**}$ and their values in the previous window.
The accuracy of $\gamma_2^{**}$ is about $10^{-9}$, but the accuracy of $\om_2^{**}$ is about $10^{-4}$, much larger. It is reasonable because the slope of $\bar \ga(\om)$ is close to zero near $(\om_2,\ga_2)$.

\begin{table}[H]
\begin{center}
\begin{tabular}{ |c|c|c|c|c| } 
\hline
window		& $\omega_2^{**}$  & $\gamma_2^{**}$ & $\De \omega_2^{**}$  & $\De \gamma_2^{**}$ \\ 
\hline
$W_1$		& 0.557003765501051 & 1.58170639609768& N/A& N/A\\
\hline
$W_2$		& 0.554773875083001 & 1.58170476081013& -2.22989041805e-3 & -1.63528755e-6\\ 
\hline
$W_3$		& 0.554837092755109 & 1.58170475989899& 6.3217672108e-5 &	-9.1114e-10 \\ 
\hline
\end{tabular}
\vspace{3mm}
\caption{$(\omega_2^{**}, \gamma_2^{**}) $ in 3 windows}\label{table2}
\end{center}
\end{table} 
\vspace{-3mm}

We could continue this zoom-in process and compute in smaller and smaller windows. However, three windows should be sufficient for illustration.

\vspace{-3mm}

\vspace{-5mm}

\subsection{The standing waves}\label{sec5.4}

In this Subsection we describe how to compute the standing waves $\phi_{\omega,\gamma}$ numerically. For one dimensional NLS \eqref{NLS} considered in this paper, because we have formula \eqref{th1.2-1} for the stability functional $J(\om,\ga)=d''(\om)$, we do not need to compute $\phi_{\omega,\gamma}$ to determine the stability of $\phi_{\omega,\gamma}$. However, when we study the same NLS in a higher dimensional setting, formula  \eqref{th1.2-1} is not available. What we can do is to compute $\phi_{\omega,\gamma}$ and 
\[
N(\omega,\gamma) = d'(\om)= \int_{\R^n}\phi_{\omega,\gamma}^2(x)\,dx,
\]
and compare it with $N(\omega+d \om,\gamma)$ to determine the sign of $d''(\om)$. Thus it is still relevant to be able to compute the standing waves $\phi_{\omega,\gamma}$ numerically.

In the 1D setting, for given $(\om,\ga)$, $\phi_{\omega,\gamma}(t)$ solves
\eqref{phi234.eq},
\EQ{\label{eq5.3.1}
\phi'' =g(\phi) = \om \phi - a_1|\phi| \phi +\ga |\phi|^{2} \phi - a_3 |\phi|^{3}\phi,
}
with $a_1,a_3 = \pm 1$, $\phi(t)>0 $ and $\lim _{t\to \pm \infty} \phi(t)=0$.
We will focus on the FDF case, $a_1=a_3=1$ and $\om,\ga>0$.  We can translate $t$ so that $\phi(0)=\max \phi$, and we have the boundary conditions
\EQ{\label{eq5.3.2}
\phi(0) = \phi_0 ,\quad \phi'(0) = 0,\quad \lim _{t\to  \infty} \phi(t)=0,
}
where $\phi_0 =\max \phi= \phi_0(\om,\ga)$ is the first positive zero of the potential function $G(t)$ given in \eqref{G234.def}
or, equivalently, the first positive zero of
\EQ{\label{phi0eq}
\om = 
\frac 25 x^{3} - \frac{\ga}2 x^{2}  +  \frac{2}3x.
}

We have tried the following three numerical methods and their combinations to compute $\phi_{\omega,\gamma}(t)$ with MATLAB.
We usually start with the time interval $0\le t \le 50$, and shift to smaller time intervals if necessary, e.g., when we take smaller $dt$ and need more computation power.

 Recall that there is no solution for $(\om,\ga)\in \Ga_\no$. For $(\om,\ga)\not\in \Ga_\no$, the computation results usually look fine. However, smaller $dt$ is required if $(\om,\ga)$ is close to $\Ga_\no$. It is because that, in this case, the solution 
of the planar dynamics \eqref{planar} with $(x(0),y(0))=(\phi_0,0)$
spends a long time in the 
neighborhood of $(\phi_0,0)$ where its velocity $(y,g(x))$ 
is extremely small. After it leaves the neighborhood, it moves rapidly toward the origin.
Hence it is a stiff computational problem when $(\om,\ga)$ is close to $\Ga_\no$.

\subsubsection{Shooting method}

In this method, we compute solutions $\phi(t) = \phi_{\om,\ga}(t)$ of the problem \eqref{eq5.3.1} with initial conditions
\EQ{%
\phi(0) = \phi_0 ,\quad \phi'(0) = 0,
}
for $t\in [0,50]$ with mesh size $\Delta t = 0.00001$, using MATLAB command \texttt{ode45}.

\begin{figure}[H]
\hspace{-10mm}
\includegraphics[width=1\textwidth, height=0.5\textwidth]{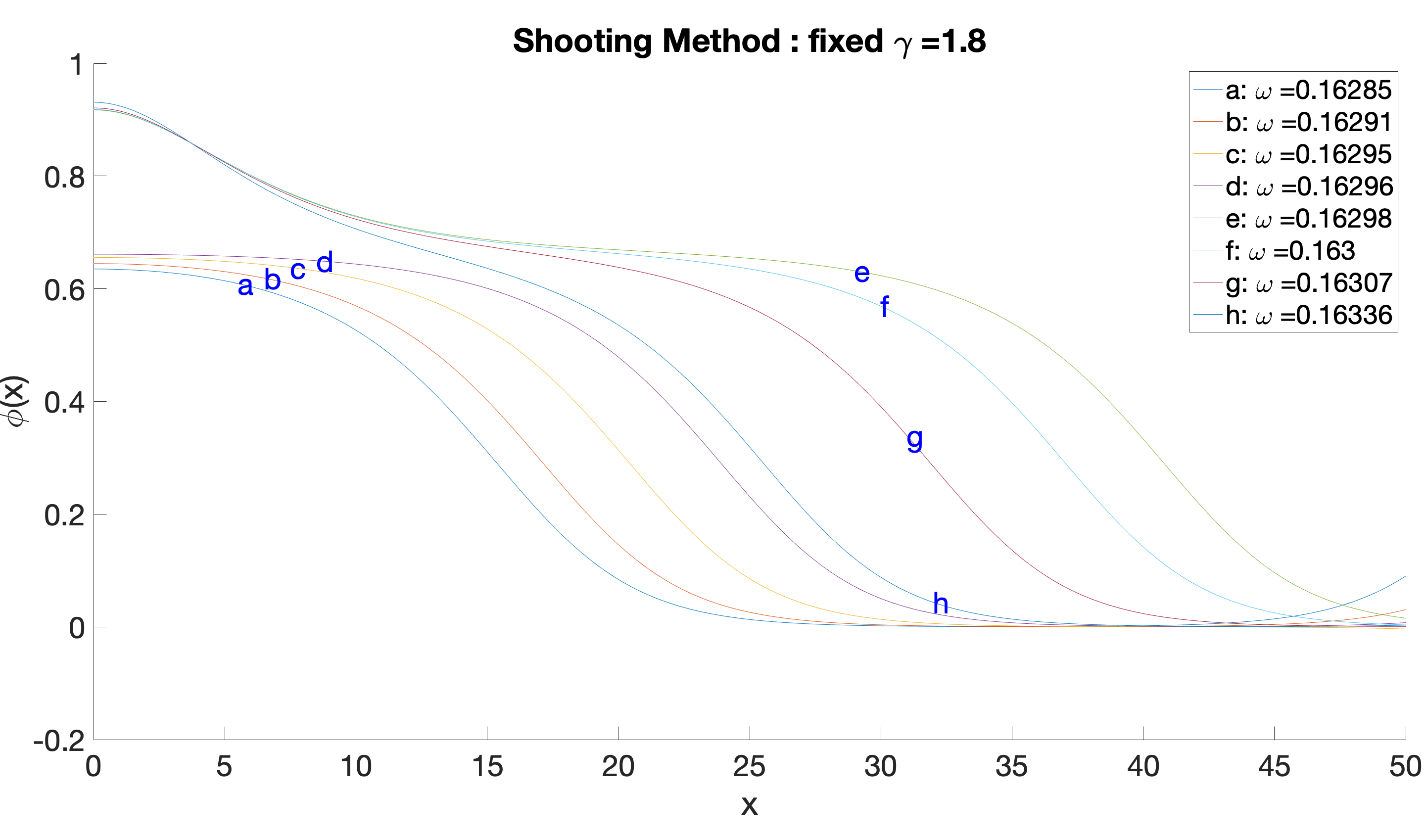}
\caption{Method 1 for $\ga=1.8$ and several $\om$'s.}
\label{Shooting}
\end{figure}

The Shooting method alone usually first gives a reasonable decaying $\phi(t)$ 
for some time. However, it then bounces back and starts oscillating, becoming a periodic solution in the long time. 
Although theoretically the solution should converge to the origin, corresponding to a homoclinic orbit passing $(\phi_0,0)$ on the phase plane, numerical errors likely perturb the solution to a nearby periodic orbit. (It looks periodic, but most likely is not, again due to numerical errors.)
Since the oscillation is due to numerical error, it should be truncated.
In Figure \ref{Shooting} we present solutions for $\ga=1.8$ and several $\om$.

Recall that $\Ga_\no$ ends at $(\omega_1,\gamma_1)\approx (0.1656, 1.7889)$ and $\ga=1.8$ is slightly larger than $\ga_1$. The values of $\om$ in Figure \ref{Shooting} are close to $\om_0= \frac {22}{135} \approx 0.162963$, which is such that $(\om_0,1.8)\in \Ga_{\no}$. The closer $\om$ to $\om_0$, the longer the solution stay near its peak $\phi_0$. As a result, the supports of the solutions in Figure \ref{Shooting} are between 25 and 50, larger than those to be found for $\ga=1.7$ in Figure \ref{Fig1.7ome}. 

\begin{figure}[H]
\centering
\includegraphics[width=0.6\textwidth, height=0.2\textwidth]{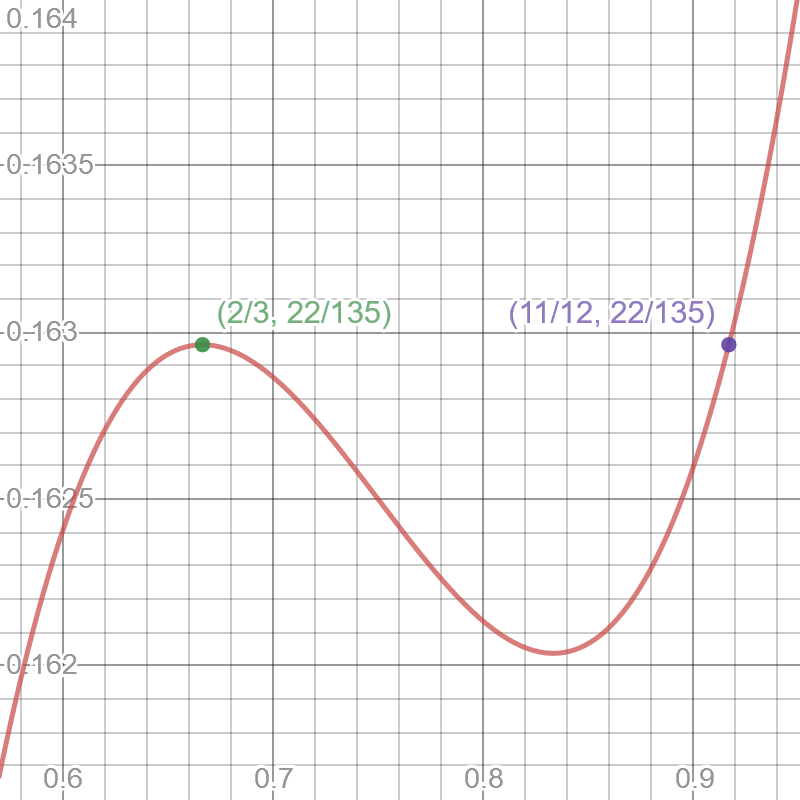}

\vspace{2mm}
\caption{Graph of the right side of \eqref{phi0eq} for the relation of $\phi_0$ and $\om$ when $\ga=1.8$.}
\label{phi0eq-gamma1d8}
\end{figure}

Recall that $\phi(0)=\phi_0$ is the first positive zero of \eqref{phi0eq} and $(\om,\ga)\in \Ga_\no$ if $\phi_0$ is a double zero of \eqref{phi0eq}. When $\ga=1.8$, $\phi_0=2/3$ is a double zero when $\om=\om_0 = \frac {22}{135}$. In this case, \eqref{phi0eq} can be factorized as
\[
\frac 25
(x-\frac 23)^2 (x-\frac {11}{12})=0,
\]
and $\frac {11}{12}$ is the third zero. See Figure \ref{phi0eq-gamma1d8} for the graph of the right side of \eqref{phi0eq} when $\ga=1.8$. Thus when $\om \to \om_0-$, $\phi_0$ is increasing and converges to $\frac 23$. 
When $\om \to \om_0+$, $\phi_0$ is decreasing and converges to $\frac {11}{12}$. There is a jump from $\frac 23$ to $\frac {11}{12}$.
This can be observed in Figure \ref{Shooting}.
In both intervals of $\om$, $(0,\om_0)$ and $(\om_0,\infty)$, $\phi_0(\om,1.8)$ is an increasing function of $\om$.

Another difficulty of the Shooting method occurs for large $\om$ and $\ga$. For example, for $(\om,\ga)=(10,5)$ and $dt=0.001$, the computation stops at $t=10.6355$ and we get the error message:
\begin{quote}
Warning: Failure at t=1.063551e+01. Unable to meet integration tolerances without reducing the step size below the smallest
value allowed (2.842171e-14) at time t. 
\end{quote}
At $t=10.6355$ the solution is small and hence has the same bahavior as the solution of the linear ODE
\[
\ddot x = \om x
\]
which has eigenvalues $\pm\sqrt \om$. 
One may imagine that a larger $\om$ corresponds to a larger positive eigenvalue and requires a smaller step size.
If we  fix $\ga = 5$, the larger the value of $\om$, the smaller the stop time. See Table \ref{table3}.

\begin{table}[H]
\begin{center}
\begin{tabular}{ |c|c|c|c|c| } 
\hline
$\om$	& 5 & 10 & 15  & 20 \\ 
\hline
stop time		& 1.330625e+01 & 1.063551e+01 & 8.375157e+00 & 7.429009e+00\\
\hline
\end{tabular}
\vspace{3mm}
\caption{stop time of shooting method for $\ga=5$ and various large $\om$}\label{table3}
\end{center}
\end{table} 
\vspace{-3mm}

A reasonable approximation is to keep the solution for $0<t<T_1$, where $T_1$ is either the first positive time that the solution has zero derivative, the time that the solution goes below zero, or the time an error occurs. One then replaces the solution by zero for $t>T_1$. We call it the \emph{shooting-cropping solution}.

\subsubsection{Picard iteration method}

In this method we refine our approximation solutions using the Picard iteration, starting from a good initial approximation. For the initial approximation we use the shooting-cropping solution from our first method. 

We now describe the Picard iteration. For the ODE $\phi'' = g(\phi)$ | given in \eqref{eq5.3.1}, 
suppose that we have a good initial guess $u_0$. 
The difference $v=\phi-u_0$ satisfies
\EQ{\label{LvFv}
Lv = F(v) =  F_0 + N(v), 
}
where
\EQN{
Lv&=v'' - g'(u_0) v , 
\\
F_0 = g(u_0)-u_0'',
\quad 
N(v)&=g(u_0+v)-g(u_0)-g'(u_0)v.
}
The source term
$F_0$ is independent of $v$ and reasonably small if $u_0$ is a good guess. The term $N(v)$ is nonlinear in $v$ and independent of $\om$.

Consider the finite difference discretization of the linear problem
\EQ{\label{linear.eq}
Lv = F, \quad \text{in}\quad (0,50) ; \quad v(0)=v(50)=0. 
}
As $\phi = u_0 + v$,
the boundary condition for $\phi$,
\EQ{\label{S5:phi.bc}
\phi(0)=\phi_0, \quad \phi(50)=0,}

is satisfied if $u_0$ also satisfies \eqref{S5:phi.bc}.
Unlike $L_0 = v'' - \om v$, problem \eqref{linear.eq} may not be invertible as $L$ may have a kernel, $L\psi=0$ for some $\psi\not=0$. However, this is non-generic, and $L$ should have no kernel if we simply perturb $u_0$ or the step-size slightly. Thus for most $u_0$,  \eqref{linear.eq} should correspond to an invertible vector equation in the finite difference method.
Let $v= G(F)$ be the solution operator of the corresponding vector equation of \eqref{linear.eq}  in the finite difference method. We can solve the solution $v$ of \eqref{linear.eq} with $F=F(v)$ by Picard iteration
\EQ{
v^{(0)}=0,\quad
v^{(k+1)} = G (F(v^{(k)})), \quad k \ge 0.
}
A revised scheme is
\EQ{
v^{(1)}=G(F_0)  ,\quad
v^{(k+1)} = v^{(1)} + G (N(v^{(k)})) , \quad k \in \N,
}
It seems more efficient as we only compute $v^{(1)}$ once.

\begin{figure}[H]
  \centering 
\includegraphics[width=0.95\textwidth]  {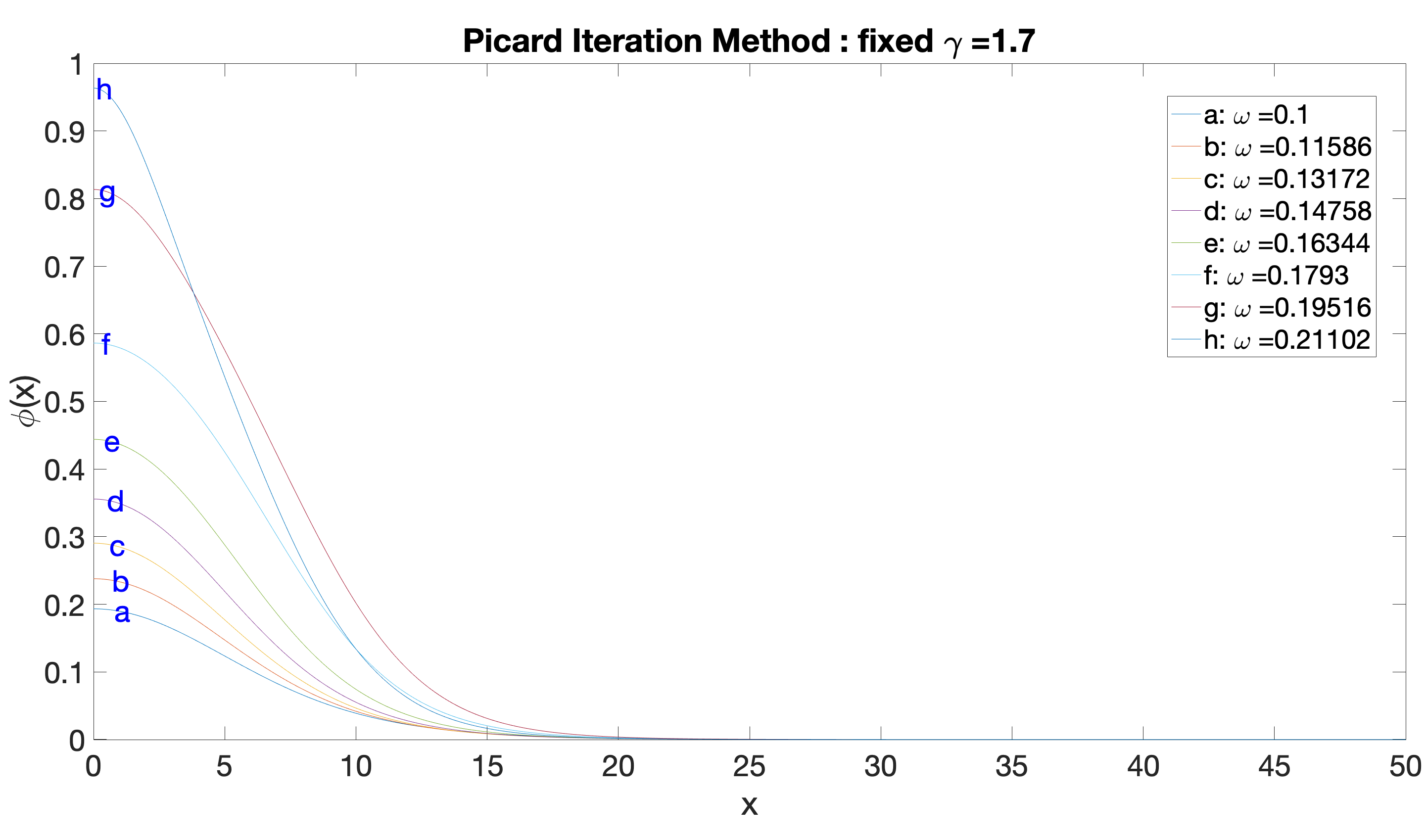}
\vspace{-2mm}
\caption{$\phi_{\omega,1.7}$ by Picard iteration method for fixed $\gamma = 1.7$ and several $\om$, with initial guess from shooting-cropping solutions and $dt = 0.1$.  }
\label{Fig1.7ome}
\end{figure}

Using this method, with dt as coarse as $dt = 0.1$, we already obtain reasonably monotonic decaying solutions which agree with the solutions obtained using $dt=0.01$, and solutions obtained using Method 3.
See Figure \ref{Fig1.7ome} for graphs by Method 2 for fixed $\gamma = 1.7$ and several $\om$. Note that $1.7<\ga_1\approx 1.7889$. Hence
$(\omega,1.7)\not \in \Ga_\no$ for all $\om$.

We also computed the solutions for $\ga=1.8$ and several choices of $\om$. Note that $\ga=1.8$ is slightly greater than $\ga_1$ and the pair $(\om,\ga)$ can be very close to $\Ga_\no$.
The result is similar to what we get using Method 3 (see Figure \ref{phisoln_gam1.8_dt00001}), hence it is skipped.

\begin{figure}[H]
\includegraphics[width=1\textwidth, height=0.5\textwidth]{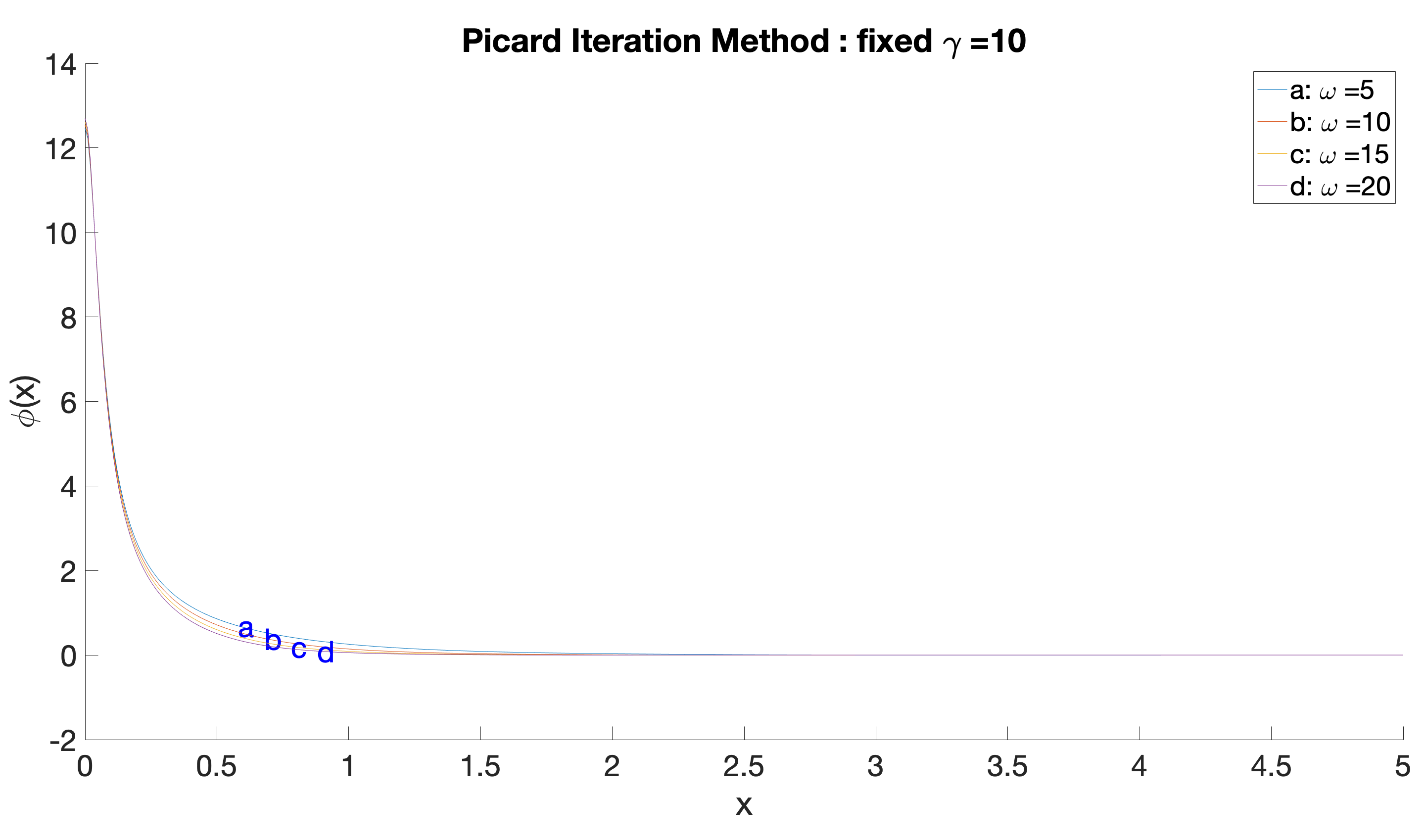}
\caption{Picard iteration method for $\ga=10$ and $\om=5,10,15,20$, with $dt = 0.01$ and the shooting-cropping solutions as initial guesses. }
\label{good_picard_example}
\end{figure}

We also computed the solutions for $\ga=10$ and several large values of $\om$. See Figure \ref{good_picard_example}. The solutions are nice decaying functions with very compact supports. The computation ending time for Figure \ref{good_picard_example} is $T=5$ for dt = 0.01 and 0.001. Our other computations with  dt = 0.01 and ending time $T=50$ give similar results.

A difficulty with the Picard iteration method is that it takes a lot of computation power for smaller $dt$. When we take $dt=0.001$, we need to shrink the time interval to at the largest $[0,20]$ to compute it in a reasonable time.

\subsubsection{MATLAB's \texttt{bvp4c} function method}

Our third method is using the MATLAB function \texttt{bvp4c} to solve the ODE $\phi'' = g(\phi)$ given in \eqref{eq5.3.1} for $0\le t \le 50$ subject to the boundary conditions \eqref{S5:phi.bc}. We use the shooting-cropping solution from Method 1 as the initial guess for the function \texttt{bvp4c}. 

\begin{figure}[H]
\includegraphics[width=0.49\textwidth, height=0.49\textwidth]{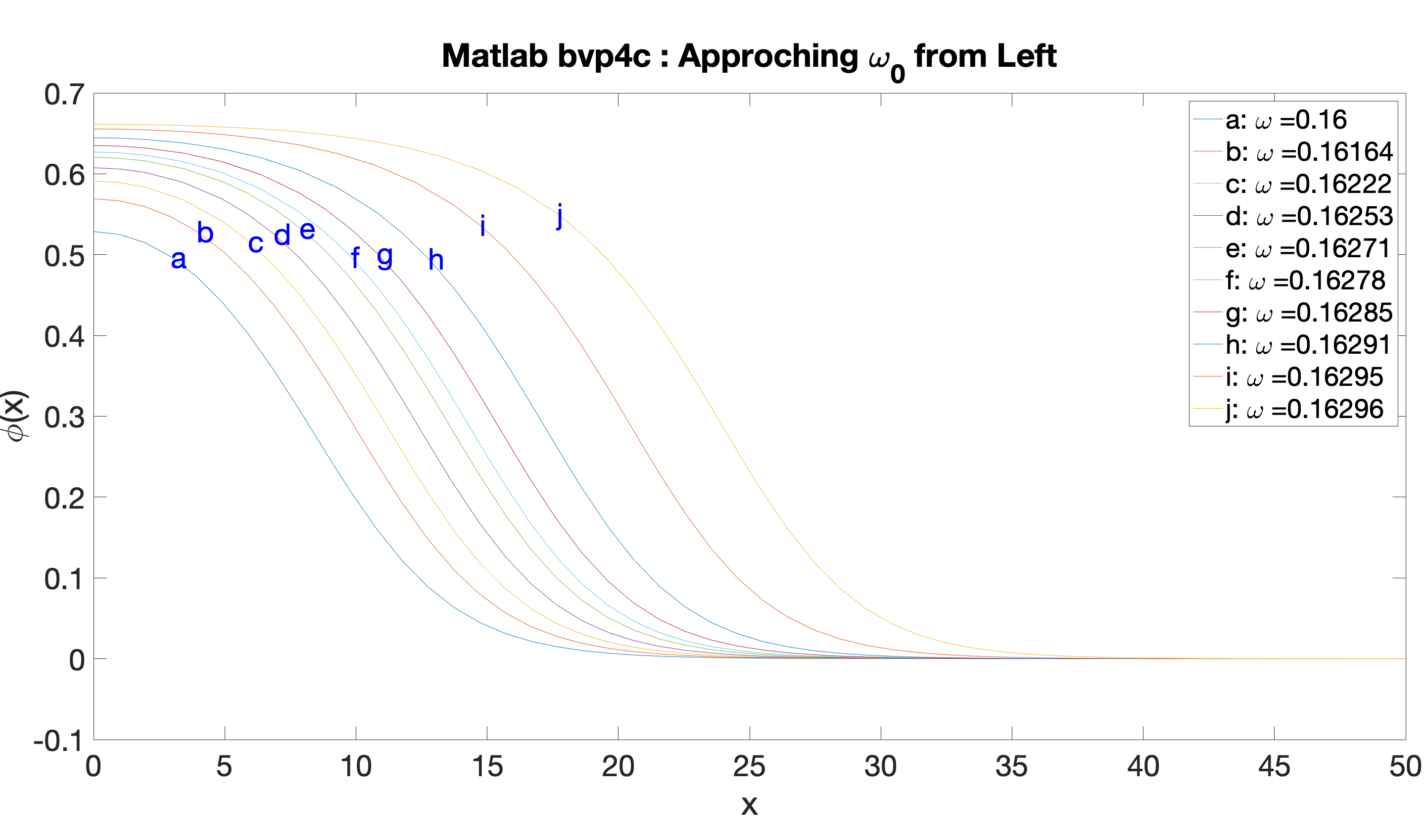}
\includegraphics[width=0.49\textwidth, height=0.49\textwidth]{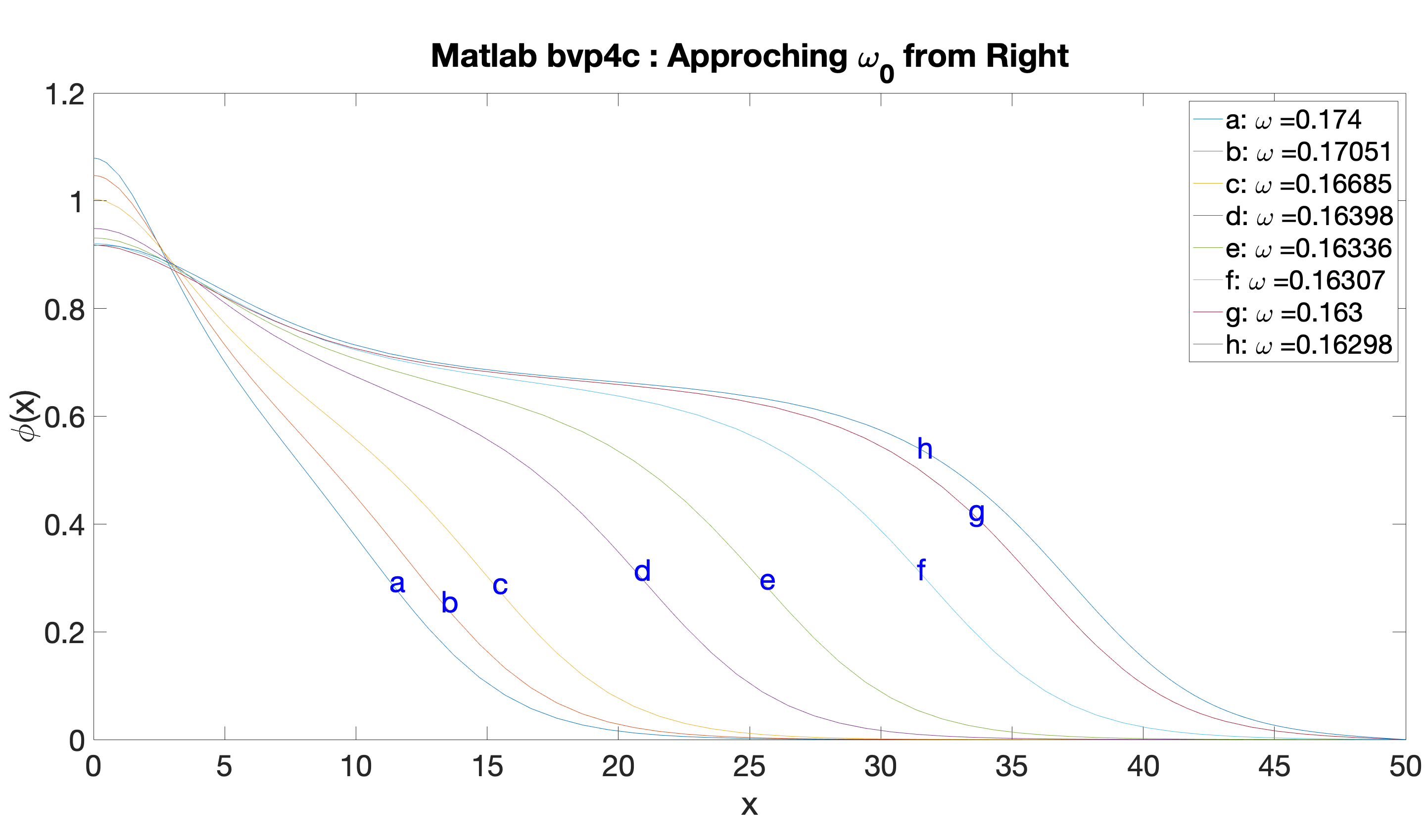}

\vspace{-2mm}
\hspace*{36mm} {\small (i) \hspace{0.46\textwidth} (ii)}
\vspace{3mm}

\caption{Method 3 for fixed $\gamma = 1.8$ for 
(i) $\omega<\omega_0$, \quad (ii) $\omega>\omega_0$.} 
\label{phisoln_gam1.8_dt00001}
\end{figure}

The numerical solutions obtained from Method 3  are usually nice looking monotonically decaying solutions. 
See Figure \ref{phisoln_gam1.8_dt00001} for a few solutions by Method 3 for $\ga=1.8$.

We also computed with Method 3 for $\ga=1.7$ and $\ga=10$. The results are similar to those obtained by Method 2, see Figures \ref{Fig1.7ome} and \ref{good_picard_example}.
Hence we skip the figures.

\section{Appendix: Explicit formulas for standing waves}
\label{Explicit-formulas}

It is well known that the solution of
\EQ{
\label{phi-p.eq}
\phi_{xx} = \om \phi -  \phi^{p} , \quad (x \in \R)
}
for $1<p<\I$ is given by 
\[
\phi_{p,\om}(x) = \om^{1/(p-1)} Q_p(\om^{1/2}x) 
\]
where $Q_p(x)= \phi_{p,1}(x)$ is given by
\[
Q_p(x) = \bke{\frac {p+1}2\, \sech^2 \bke{\frac {p-1}2 x}}^{\frac 1{p-1}}
\]

In the following we consider explicit solutions for double power nonlinearities.
Let
\EQ{
v(x) = (\ell+\ch^2 x)^{-\beta}
}
where $\ell\in(-1,\I)$ and $\beta>0$. Using $\sh^2 = \ch^2-1$, one gets
$v_x(x) =-2\beta (\ell+\ch^2 x)^{-\beta-1}\ch x \,\sh x$ and
\[
v_{xx} = Av - B v^{1+1/\beta} - C  v^{1+2/\beta},
\]
where
\[
A= 4\beta^2, \quad B=2\beta(2\beta+1)(2\ell+1),\quad C=-4\beta(\beta+1)\ell(\ell+1).
\]
To get double power nonlinearity we require $BC \not =0$, i.e., $\ell \not = -1/2, 0$.
There are three cases:
\begin{enumerate}
\item $\ell\in I_1=(-1,-1/2)$: We have $B<0$, $C>0$, (defocusing-focusing, DF)

\item $\ell\in I_2=(-1/2,0)$: We have $B>0$, $C>0$, (focusing-focusing, FF)

\item $\ell\in I_3=(0,\I)$: We have $B>0$, $C<0$, (focusing-defocusing, FD)
\end{enumerate}

\begin{remark}
The borderline cases $\ell=-1/2$ and  $\ell=0$ correspond to NLS with a focusing single power nolinearity, and suggest that the
profile of $v$ near $\ell=-1/2$ and  $\ell=0$ are given by $Q_{1+1/\beta}$ and $Q_{1+2/\beta}$, respectively.
\end{remark}

We now assume $\ell\not = -1/2,0$ and
denote $a_1=\sgn B$ and $a_2 = \sgn C$.
A suitable rescaling
\[
\phi(x) = k v(\la x), \quad k=|C/B|^\beta, \quad \la = \frac{\sqrt{|C|}}{|B|}, \quad
\om = 4 \beta^2\la^2
\]
gives a solution $\phi(x)$ of
\EQ{
\label{35phi.eq}
\phi_{xx} = \om \phi -  a_1\phi^{1+1/\beta}  -a_2  \phi^{1+2/\beta}.
}
Explicitly, 
\[
\phi_\om(x) = k\bke{\ell+  \ch^2 (\frac {\sqrt {\om}}{2 \beta}x) }^{-\beta},
\quad k=\abs{ \frac {2(\beta+1)\ell(\ell+1)}{ (2\beta+1)(2\ell+1)}}^\beta,
\]
and
\EQ{
\label{om*.def}
\om = \om^*\frac{4|\ell|(\ell+1)}{(2\ell+1)^2}, \quad
\om^*=\frac{\beta(\beta+1)}{(2\beta+1)^2}.
}
\begin{remark}
Explicit solutions for \eqref{35phi.eq} are well known. For example \cite[(3.1)]{Oh95} 
has an %
equivalent formula for solutions of \eqref{35phi.eq} with $a_1,a_2 \in \R$:
\[
\phi_\om(x) = \bke{\frac {\om}{A+\sqrt{A^2+B\om}\,\ch(\beta^{-1}\sqrt {\om}x)}}^{\beta},
\quad A=\frac {a_1}{2+1/\beta}, \quad B=\frac {a_2}{1+1/\beta}.
\]
The main advantage of our form is the single parameter $\ell$ for all three cases: DF, FF and FD.
\end{remark}

The original parameter $\ell$ can be solved in each interval $I_j$ in terms of
$\om$. Thus we can use $\om$ as the parameter in each interval.  Indeed,
from \eqref{om*.def}, we have
\[
\om_{-1+}=0, \quad \om_{-1/2-} = \I = \om_{-1/2+}, \quad \om_{0-} = 0 = \om_{0+}, 
\quad \om(\I-)=\om^*,
\]
and
\[
\frac {d\om}{d\ell} = \frac {4 \om^* \sgn \ell}{(2\ell+1)^3}
\]
which is positive for $\ell \in (-1,-1/2)\cup (0,\I)$, and negative for
$\ell \in (-1/2,0)$.

\begin{figure}[H]
\hspace{10mm}
\includegraphics[width=\textwidth, height=0.4\textwidth]{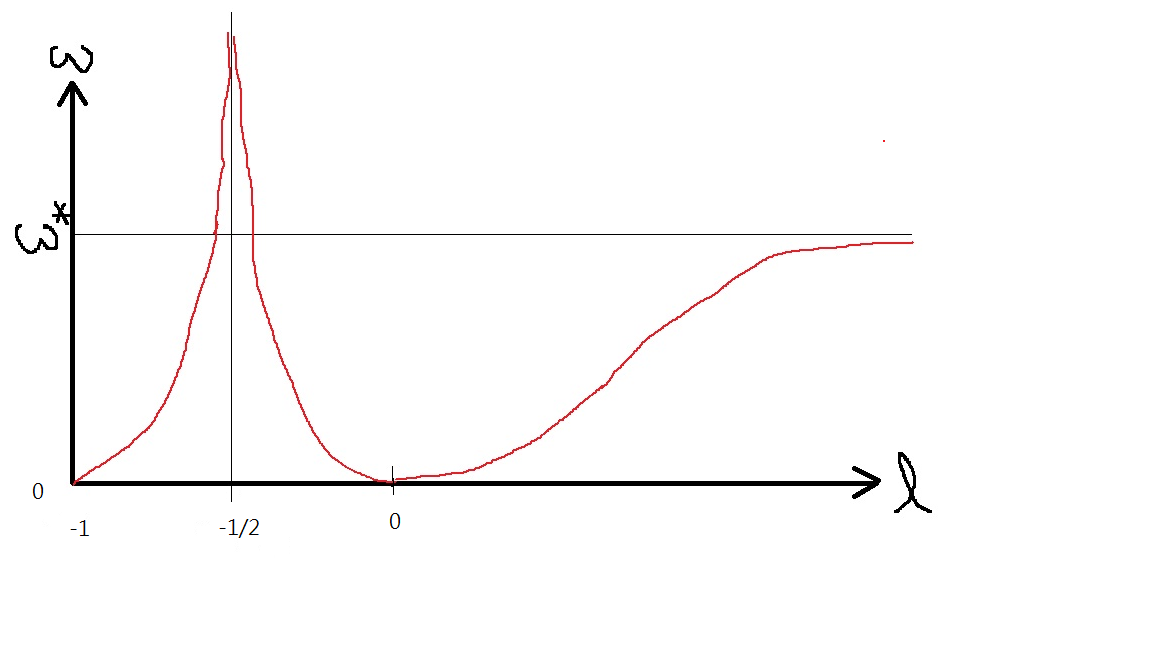}
\vspace{-12mm}
\caption{$\omega$ as a function of $\ell$}
\end{figure}

\begin{remark}\label{rem1.6}
Consider the limit $\ell= -1+\e$ with $\e \to 0+$. We have $\om \sim C\e$ and
$\phi_\om(x) \sim \phi_0(x):=c_1(c_2+x^2)^{-\beta}$ for $|x|< \e^{-1/2}$ and $\phi_\om(x) \sim c_3 \e^\beta e^{-c_4\sqrt {\e}|x|}$ for $|x|> \e^{-1/2}$.
It can be understood as the competition between $f(\phi)$ and $\om \phi$ for $\om \ll 1$. When $|f(\phi)| > \om \phi$,  i.e., when $|x|< \e^{-1/2}$, 
 $\phi_\om$ is approximated by
$\phi_0$, the solution of $\phi''+f(\phi)=0$. For $|x|> \e^{-1/2}$, $\om \phi $ is larger than $f(\phi)$ and $\phi_\om$  is approximated the solution of $\phi''=\om \phi$. Similarly, we can consider the limits $\ell \to 0_-$ and $\ell \to 0_+$.
\end{remark}

\section*{Acknowledgments}
We warmly thank Vianney Combet for helpful discussions and continued interests in this work. We thank Stefan Le Coz for the reference \cite{Fukaya-Hayashi}. We also thank the referee for very valuable suggestions. The work of Tsai was partially supported by NSERC grant RGPIN-2018-04137.

\addcontentsline{toc}{section}{\protect\numberline{}{References}}
\bibliographystyle{abbrv}\bibliography{234nls2021}

\end{document}